\documentstyle[11pt]{article}
\newfam\msbfam
\font\tenmsb=msbm10    \textfont\msbfam=\tenmsb \font\sevenmsb=msbm7
\scriptfont\msbfam=\sevenmsb \font\fivemsb=msbm5
\scriptscriptfont\msbfam=\fivemsb
\def\Bbb{\fam\msbfam \tenmsb}

\def\rr{{\Bbb R}}
\def\rz{{{\rr}^n}}
\def\zz{{\Bbb Z}}

\def\fz{\infty}
\def\az{\alpha}

\def\ez{\epsilon}
\def\bz{\beta}

\def\tz{\theta}
\def\sz{\sigma}
\def\vz{\varphi}
\def\dz{\delta}
\def\gz{\gamma}

\def\lz{\lambda}
\def\oz{\Omega}

\def\wt{\widetilde}
\def\wz{\omega}
\def\l{\left}
\def\r{\right}

\def\dsum{\displaystyle\sum}
\def\dprod{\displaystyle\prod}
\def\dint{\displaystyle\int}
\def\dfrac{\displaystyle\frac}
\def\dsup{\displaystyle\sup}

\def\dinf{\displaystyle\inf}

\newtheorem{thm}{\hskip\parindent Theorem}

\newtheorem{lem}{\hskip\parindent Lemma}

\newtheorem{prop}{\hskip\parindent Proposition}

\newtheorem{cor}{\hskip\parindent Corollary}

\begin{document}

\baselineskip=15pt
\renewcommand{\arraystretch}{2}
\arraycolsep=1.2pt

\title{ Extrapolation from $A_\fz^{\rho,\fz}$, vector-valued inequalities and applications in the  Schr\"odinger settings
 \footnotetext{ \hspace{-0.65
cm} 2000 Mathematics Subject  Classification: 42B25, 42B20.\\
The  research was supported  by the NNSF (10971002) of China.\\}}

\author{ Lin Tang }
\date{}
\maketitle

{\bf Abstract}\quad   In this paper, we generalize the $A_\fz$ extrapolation
theorem in \cite{cmp} and the $A_p$ extrapolation theorem of Rubio
de Francia to Schr\"odinger settings. In addition, we also establish
the weighted vector-valued inequalities for Schr\"odinger type maximal operators by using weights belonging to $ A_p^{\rho,\tz}$ which includes $A_p$. As their applications, we establish the weighted vector-valued inequalities for some
Sch\"odinger type operators and pseudo-differential operators.
\bigskip

\begin{center}{\bf 1. Introduction }\end{center}

In this paper, we consider the  Sch\"odinger differential operator
$L=-\Delta+V(x)\ {\rm on}\ \rz,$ $\ n\ge 3,$ where  $V(x)$
is a  nonnegative potential satisfying  certain reverse H\"older
class.

 We say a nonnegative locally $L^q$
integral function $V(x)$ on $\rr^n$ is said to belong to $B_q(1<q\le
\fz)$ if there exists $C>0$ such that the reverse H\"older
inequality
$$\l(\dfrac 1{|B(x,r)|}\dint_{B(x,r)} V^q(y)dy\r)^{1/q}\le C\l(\dfrac
1{|B(x,r)|}\dint_{B(x,r)} V(y)dy\r)$$ holds for every $x\in\rz$ and
$0<r<\fz$, where $B(x,r)$ denotes the ball centered at $x$ with
radius $r$. In particular, if $V$ is a nonnegative polynomial, then
$V\in B_\fz$.  Throughout this paper, we always assume that
$0\not\equiv V\in B_n/2$.

The study of schr\"odinger operator $L=-\triangle+V$ recently
attracted much attention; see \cite{bhs1,bhs2,dz,dz1,s1,z}.
In particular, it should be pointed out that Shen \cite{s1} proved
the Schr\"odinger type operators, such as
$\nabla(-\Delta+V)^{-1}\nabla$, $\nabla(-\Delta+V)^{-1/2}$,
$(-\Delta+V)^{-1/2}\nabla$  with $V\in B_n$, $(-\Delta+V)^{i\gz}$  with $\gz\in\rr$ and $V\in B_{n/2}$,  are standard
Calder\'on-Zygmund operators.

 Recently, Bongioanni, etc,  \cite{bhs1} proved
$L^p(\rz)(1<p<\fz)$ boundedness for commutators of Riesz transforms
associated with Schr\"odinger operator with
$BMO_\tz(\rho)$ functions which include the class $BMO$ function,
and in \cite{ bhs2} established the weighted boundedness
for
Riesz transforms, fractional integrals  and Littlewood-Paley functions
 associated with Schr\"odinger operator with
weight $A_p^{\rho,\tz}$ class which includes the Muckenhoupt weight class.
Very recently, the author \cite{t1, t2} established
the weighted norm inequalities for some Schr\"odinger type operators,
 which include commutators of Riesz transforms, fractional
integrals and Littlewood-paley operators.

On the other hand, extrapolation for weights plays an important role
in Harmonic analysis. In particulary, Rubio de Francia \cite{ru}
proved the $A_p$ extrapolation theorem: If the operator is bounded
on $L^{p_0}(\wz)$ for some $p_0$, $1<p_0<\fz$, and every $\wz\in
A_{p_0}$, then for every $p$, $1<p<\fz$, $T$ is bounded on
$L^p(\wz), \wz\in A_p$ (see also \cite{d, g}). Recently,
Cruz-Uribe, etc, in \cite{cmp} extended this theorem from $A_p$
weights to $A_\fz$ weights, to pairs of operators, and to the range
$0<p<\fz$ in the context of Muckenhoupt bases.

In this paper, we generalize the $A_\fz$ extrapolation
theorem in \cite{cmp} and the $A_p$ extrapolation theorem of Rubio
de Francia to Schr\"odinger settings and give some applications.

The paper is organized as follows. In Section 2,  we give
factorization of $A_p^{\rho,\fz}$, and establish the
weighted vector-valued inequalities for Schr\"odinger type maximal operators, these
results play a crucial role in this paper.  In Section 3, we
obtain extrapolation theorems from $A_\fz^{\rho,\fz}$ and $A_p^{\rho,\fz}$. Finally, we establish the
weighted vector-valued inequalities for some Schr\"odinger type
operators and pseudo-differential operators in section 4.

Throughout this paper, we let $C$ denote  constants that are
independent of the main parameters involved but whose value may
differ from line to line. By $A\sim B$, we mean that there exists a
constant $C>1$ such that $1/C\le A/B\le C$.

\bigskip

\begin{center} {\bf 2. Factorization and vector-valued inequalities } \end{center}
In this section, we give the factorization of $A_p^{\rho,\fz}$ and weighted
vector-valued inequalities for  Schr\"odinger type maximal operators.

We first recall some notation.  Given $B=B(x,r)$
and $\lz>0$, we will write $\lz B$ for the $\lz$-dilate ball,
which is the ball with the same center $x$ and with radius $\lz
r$. Similarly, $Q(x,r)$ denotes the cube centered at $x$ with the
sidelength $r$ (here and below only cubes
with sides parallel to the coordinate axes are considered), and $\lz Q(x,r)=Q(x,\lz r)$.
 Let $f=\{f_k\}_1^\fz$ is a sequence of locally integral functions $\rz$,
 $|f(x)|_r=(\sum_{k=1}^\fz|f_k(x)|^r)^{1/r}$, and $|Tf(x)|_r=(\sum_{k=1}^\fz|Tf_k(x)|^r)^{1/r}$.

 The function $m_V(x)$ is defined by
$$\rho(x)=\dfrac 1{m_V(x)}=\dsup_{r>0}\l\{r:\ \dfrac
{1}{r^{n-2}}\dint_{B(x,r)}V(y)dy\le 1\r\}.$$Obviously,
$0<m_V(x)<\fz$ if $V\not=0$. In particular, $m_V(x)=1$ with $V=1$
and $m_V(x)\sim (1+|x|)$ with $V=|x|^2$.

\begin{lem}\label{l2.1.}\hspace{-0.1cm}{\rm\bf 2.1(\cite{s1}).}\quad
There exists $l_0>0$ and $C_0>1$such that $$\dfrac
1{C_0}\l(1+|x-y|m_V(x)\r)^{-l_0}\le \dfrac{m_V(x)}{m_V(y)}\le
C_0\l(1+|x-y|m_V(x)\r)^{l_0/(l_0+1)}.
$$ In particular, $m_V(x)\sim m_V(y)$ if $|x-y|<C/m_V(x)$.
\end{lem}
In this paper, we write $\Psi_\tz(B)=(1+r/\rho(x_0))^\tz$, where
 $\tz> 0$, $x_0$ and $r$ denotes the center and radius of $B$ respectively.

A weight will always mean a positive function which is locally
integrable. As \cite{bhs2}, we say that a weight $\wz$ belongs to the class
$A_p^{\rho,\tz}$ for $1<p<\fz$, if there is a constant $C$ such that for
all balls
  $B$
$$\l(\dfrac 1{\Psi_\tz(B)|B|}\dint_B\wz(y)\,dy\r)
\l(\dfrac 1{\Psi_\tz(B)|B|}\dint_B\wz^{-\frac 1{p
-1}}(y)\,dy\r)^{p-1}\le C.$$
 We also
say that a  nonnegative function $\wz$ satisfies the $A_1^{\rho,\tz}$
condition if there exists a constant $C$ such that
$$M_{V,\tz}(\wz)(x)\le C \wz(x), \ a.e.\ x\in\rz.$$
 where
$$M_{V,\tz}f(x)=\dsup_{x\in B}\dfrac 1{\Psi_\tz(B)|B|}\dint_B|f(y)|\,dy.$$
When $V=0$, we
denote $M_{0}f(x)$ by $Mf(x)$( the standard Hardy-Littlewood
maximal function). It is easy to see that $|f(x)|\le M_{V,\tz} f(x)\le
Mf(x)$ for $a.e.\ x\in\rz$ and any $\tz\ge 0$.

Since $\Psi_\tz(B)\ge 1$ with $\tz\ge 0$, then
$A_p\subset A_p^{\rho,\tz}$ for $1\le p<\fz$, where $A_p$ denotes
the classical Muckenhoupt weights; see \cite{gr} and \cite{m}. We
will see that $A_p\subset\subset A_p^{\rho,\tz}$ for $1\le p<\fz$ in
some cases.  In fact, let $\tz>0$ and $0\le\gz\le\tz$, it is easy
to check that $\wz(x)=(1+|x|)^{-(n+\gz)}\not\in A_\fz$ and
$\wz(x)dx$ is not a doubling measure, but
$\wz(x)=(1+|x|)^{-(n+\gz)}\in A_1^{\rho,\tz}$ provided that  $V=1$ and
$\Psi_\tz (B(x_0,r))=(1+r)^\tz$.

We remark that  balls can be replaced by cubes  in definition of
$A_p^{\rho,\tz}$ and $M_{V,\tz}$, since
$\Psi(B)\le \Psi(2B)\le 2^\tz \Psi(B).$

Next we give the weighted boundedness of $M_{V,\tz}$.
\begin{lem}\label{l2.2.}\hspace{-0.1cm}{\rm\bf 2.2(\cite{t}).}\quad
Let $1<p<\fz$, $p'=p/(p-1)$ and assume that $\wz\in A_p^{\rho,\tz}$.
 There
exists a constant $C>0$ such that
$$\|M_{V,p'\tz}f\|_{L^p(\wz)}\le C\|f\|_{L^p(\wz)}.$$
\end{lem}

Similar to the classical Muckenhoupt weights(see \cite{j,gr,s}), we give some
properties for weight class $A_p^{\rho,\tz}$ for $p\ge 1$.
\begin{prop}\label{p2.1.}\hspace{-0.1cm}{\rm\bf 2.1.}\quad
Let $\wz\in A^{\rho,\fz}_p=\bigcup_{\tz\ge
0}A_p^{\rho,\tz}$ for $p\ge 1$. Then
\begin{enumerate}
\item[(i)]If $ 1\le p_1<p_2<\fz$, then $A_{p_1}^{\rho,\tz}\subset
A_{p_2}^{\rho,\tz}$. \item[(ii)] $\wz\in A_p^{\rho,\tz}$ if and only
if $\wz^{-\frac 1{p-1}}\in A_{p'}^{\rho,\tz}$, where $1/p+1/p'=1.$
\item[(iii)] If $\wz\in A_p^{\rho,\fz},\ 1<p<\fz$, then there exists
$\ez>0$ such that $\wz\in A_{p-\ez}^{\rho,\fz}.$
\item[(vi)] Let $f\in L_{loc}(\rz)$, $0<\dz<1$, then
$(M_{V,\tz})^\dz\in A_1^{\rho,\tz}$.
\item[(v)] Let $1<p<\fz$, then $\wz\in A_p^{\rho,\fz}$ if and only
if $\wz=\wz_1\wz_2^{1-p}$, where $\wz_1,\wz_2\in A_1^{\rho,\fz}$.
\end{enumerate}
\end{prop}
{\it Proof.}\quad (i) and (ii) are obvious by the definition of
$A_p^{\rho,\tz}$. (iii) is proved in \cite{bhs2}. In fact, from
Lemma 5 in \cite{bhs2}, we know that if $\wz\in A_p^{\rho,\tz}$,
then $\wz\in A_{p_0}^{\rho,\tz_0}$, where $p_0=1+\frac
{p-1}{1-\dz}<p$ with $0<\dz<1$($\dz$ is a constant depending only on
the $A_p^{\rho,loc}$ constant of $\wz$, see \cite{bhs2}) and $\tz_0=\frac {\tz p+\eta
(p-1)}{p_0}$ with $\eta=\tz
p+(\tz+n)\frac{pl_0}{l_0+1}+(l_0+1)\frac{n\dz}{1+\dz}$. We now prove
(vi). It will suffice to show that there exists a constant $C$ such
that for every $f$, every cube $Q$ and almost every $x\in Q$,
$$\dfrac1{\Psi_\tz(Q)|Q|}\dint_QM_{V,\tz}f(y)^\dz dy\le
CM_{V,\tz}f(x)^\dz.$$ Fix $Q$ and decompose $f$ as $f=f_1+f_2$,
where $f_1=f\chi_{2Q}$. Then $M_{V,\tz}f(x)\le M_{V,\tz}f_1(x)+
M_{V,\tz}f_2(x)$, and so for $0\le\dz<1$,
$$M_{V,\tz}f(x)^\dz\le M_{V,\tz}f_1(x)^\dz+ M_{V,\tz}f_2(x)^\dz.$$
Since $M_{V,\tz}$ is weak (1,1), by Kolmogorev's inequality( see \cite{rr})
$$\begin{array}{cl}
\dfrac1{\Psi_\tz(Q)|Q|}\dint_Q(M_{V,\tz}f_1)^\dz(y)dy&\le \dfrac
C{\Psi_\tz(Q)|Q|}|Q|^{1-\dz}\|f_1\|_1^\dz\\
&\le C\l(\dfrac 1{\Psi_\tz(Q)|Q|}\dint_{2Q}|f(y)|dy\r)^\dz\\
&\le C\l(\dfrac 1{\Psi_\tz(2Q)|2Q|}\dint_{2Q}|f(y)|dy\r)^\dz \\
&\le C M_{V,\tz}f(x)^\dz.
\end{array}$$
To estimate $M_{V,\tz}f_2$, note that let $Q'$ is a cube such that
$x\in Q'$, if  $Q'\bigcap(\rz\setminus (2Q))\not=\emptyset$, then $Q\subset
4nQ'$. Hence, for any $z\in Q$
$$\dfrac1{\Psi_\tz(Q')|Q'|}\dint_{Q'}|f_2(y)|dy\le \dfrac C{\Psi_\tz(4nQ')|4nQ'|}\dint_{4nQ'}|f_2(y)|dy\le
CM_{V,\tz}(z).$$ So $M_{V,\tz}(y)\le CM_{V,\tz}(x)$ for any $y\in
Q$. Thus
$$\dfrac1{\Psi_\tz(Q')|Q'|}\dint_{Q'}M_{V,\tz}f_2(y)^\dz dy\le C
M_{V,\tz}f(x)^\dz.$$ It remains to prove (v). We first assume
$\wz_1\in A_1^{\rho,\tz_1}$ and  $\wz_2\in A_1^{\rho,\tz_2}$. Since
$$\l(\dfrac
1{\Psi_{\tz_1}(Q)|Q|}\dint_{Q}\wz_1(y)dy\r)\l(\dinf_{Q}\wz_1(y)\r)^{-1}\le
C_1,$$
$$\l(\dfrac
1{\Psi_{\tz_2}(Q)|Q|}\dint_{Q}\wz_2(y)dy\r)\l(\dinf_{Q}\wz_2(y)\r)^{-1}\le
C_2,$$ moreover
$$\begin{array}{cl}
\dfrac1{\Psi_{\tz}(Q)|Q|}\dint_{Q}\wz(y)dy&=\dfrac1{\Psi_{\tz}(Q)|Q|}\dint_{Q}\wz_1(y)\wz_2^{1-p}(y)dy\\
&\le\l(\dfrac
1{\Psi_{\tz}(Q)|Q|}\dint_{Q}\wz_1(y)dy\r)\l(\dinf_{Q}\wz_2(y)\r)^{1-p},\end{array}$$
$$\begin{array}{cl}
\l(\dfrac1{\Psi_{\tz}(Q)|Q|}\dint_{Q}\wz^{-\frac1{p-1}}(y)dy\r)^{p-1}&=
\l(\dfrac1{\Psi_{\tz}(Q)|Q|}\dint_{Q}\wz_1^{-\frac1{p-1}}(y)\wz_2(y)dy\r)^{p-1}\\
&\le\l(\dfrac
1{\Psi_{\tz}(Q)|Q|}\dint_{Q}\wz_2(y)dy\r)^{p-1}\l(\dinf_{Q}\wz_1(y)\r)^{-1}.\end{array}$$
From these inequalities above and choosing
$\tz=\max\{\tz_1,\tz_2\}$, then
$$\l(\dfrac1{\Psi_{\tz}(Q)|Q|}\dint_{Q}\wz(y)dy\r)\l(\dfrac1{\Psi_{\tz}(Q)|Q|}\dint_{Q}\wz^{-\frac1{p-1}}(y)dy\r)^{p-1}\le
C_1C_2^{p-1}.$$ To prove the converse, we consider first $p\ge2$,
let $\wz\in A_p^{\rho,\tz}$, and define $T$ by
$$Tf=[\wz^{-1/p}M_{V,p\tz}(f^{p/p'}\wz^{1/p})]^{p'/p}+\wz^{1/p}M_{V,p\tz}(f\wz^{-1/p}).$$
Because $\wz^{-p'/p}\in A_{p'}^{\rho,\tz}$, then $T$ is bounded on
$L^p$ by Lemma 2.2, that is,
$$\|Tf\|_{L^p}\le A\|f\|_{L^p},$$ for some $A>0$. Also,
since $p\ge 2, p/p'\ge1$, and  Minkowski's inequality gives
$T(f_1+f_2)\le Tf_1+Tf_2$. Fix now a nonnegative $f$
with$\|f\|_{L^p}=1$ and write
$$\eta=\dsum_{k=1}^\fz(2A)^{-k}T^k(f),$$
where $T^k(f)=T(T^{k-1}(f))$. Then $\|\eta\|_{L^p}\le 1$.
Furthermore, since $T$ is positivity-preserving and subadditive, we
have the pointwise inequality
$$T\eta\le
\dsum_{k=1}^\fz(2A)^{-k}T^{k+1}(f)=\dsum_{k=2}^\fz(2A)^{1-k}T^k(f)\le
2A\eta.$$ Thus, if $\wz_1=\wz^{1/p}\eta^{p/p'}$, then
$$M_{V,p\tz}(\wz_1)\le (T(\eta))^{p/p'}\wz^{1/p}\le (2A\eta )^{p/p'}\wz^{1/p}=(2A )^{p/p'}\wz_1$$
and $\wz\in A_1^{\rho,p\tz}$. Similarly, if $\wz_2=\wz^{-1/p}\eta$,
then $M_{V,p\tz}(\wz_1)\le 2A\wz_2$, so $\wz_2\in A_1^{\rho,p\tz}$.
Moreover, $$\wz=
\wz_1\wz_2^{1-p}=\wz^{1/p}\eta^{p/p'}(\wz^{-1/p}\eta)^{1-p},$$ since
$p/p'=p-1$, finishing the proof or $p\ge 2$.

The case $p\le 2$ is similar. In fact, let $\wz\in A_p^{\rho,\tz}$,
then $\wz^{-p'/p}\in A_{p'}^{\rho,\tz}$,  and define $T$ by
$$Tf=[\wz^{1/p}M_{V,p'\tz}(f^{p'/p}\wz^{-1/p})]^{p/p'}+\wz^{-1/p}M_{V,p'\tz}(f\wz^{1/p}).$$
then $T$ is bounded on $L^p$ by Lemma 2.2, that is,
$$\|Tf\|_{L^{p'}}\le B\|f\|_{L^{p'}},$$ for some $A>0$ Also,
since $p\le 2, p'/p\ge1$, and  Minkowski's inequality gives
$T(f_1+f_2)\le Tf_1+Tf_2$. Fix now a nonnegative $f$
with$\|f\|_{L^{p'}}=1$ and write
$$\eta=\dsum_{k=1}^\fz(2B)^{-k}T^k(f),$$
where $T^k(f)=T(T^{k-1}(f))$. Then $\|\eta\|_{L^{p'}}\le 1$.
Furthermore, since $T$ is positivity-preserving and subadditive, we
have the pointwise inequality
$$T\eta\le
\dsum_{k=1}^\fz(2B)^{-k}T^{k+1}(f)=\dsum_{k=2}^\fz(2B)^{1-k}T^k(f)\le
2B\eta.$$ Thus, if $\wz_1=\wz^{-1/p}\eta^{p'/p}$, then
$$M_{V,p\tz}(\wz_1)\le (T(\eta))^{p'/p}\wz^{-1/p}\le (2B\eta )^{p'/p}\wz^{1/p}=(2B )^{p'/p}\wz_1$$
and $\wz\in A_1^{\rho,p'\tz}$. Similarly, if $\wz_2=\wz^{1/p}\eta$,
then $M_{V,p'\tz}(\wz_1)\le 2B\wz_2$, so $\wz_2\in
A_1^{\rho,p'\tz}$. Moreover, $$\wz=
\wz_2\wz_1^{1-p}=\wz^{1/p}\eta(\wz^{-1/p}\eta^{p'/p})^{1-p},$$ since
$p/p'=p-1$, finishing the proof or $p\le 2$. The proof is complete.\hfill$\Box$

C. Fefferman and E. Stein \cite{fs} obtained the vector-valued inequalities of Hardy-Littlewood maximal operators. Later, K. Andersen and R. John \cite{aj} generalized the Fefferman-Stein vector-valued inequalities to $A_p$ weights case. We next give some weighted vector-valued inequalities of maximal operators $M_{V,\eta}$ by new weights. The following interpolation results will be required. Let ${\cal S}$
denote the linear space of sequence $f=\{f_k\}$ of the form:
$f_k(x)$ is a simple function on $\rz$ and $f_k(x)\equiv 0$ for all
sufficient large $k$. ${\cal S}$ is dense in $L_\wz^p(l^r)$, $1\le
p, r<\fz$; see \cite{bp}.
\begin{lem}\label{l2.3.}\hspace{-0.1cm}{\rm\bf 2.3(\cite{aj}).}\quad
Let $\wz\ge 0$ be locally integral on $\rz$, $1<r<\fz$, $1\le p_i\le
q_i<\fz$ and suppose $T$ is a sublinear operator defined on ${\cal
S}$ satisfying
$$\wz(\{x\in\rz: |Tf(x)|_r>\az\})\le
M_i^{q_i}\az^{-q_i}\l(\int_\rz|f(x)|_r^{p_i}\wz(x)dx\r)^{q_i/p_i}$$
for $i=0,1$ and $f\in {\cal S}$. Then $T$ extends uniquely to a
sublinear operator on $L_\wz^p(l^r)$ and there is a constant $M_\tz$
such that
$$\l(\int_\rz|Tf(x)|_r^q\wz(x)dx\r)^{1/q}\le M_\tz\l(\int_\rz|f(x)|_r^{p}\wz(x)dx\r)^{1/p}$$
where $(1/p,1/q)=(1-\tz)(1/p_0,1/q_0)+\tz(1/p_1.1/q_1),\quad
0<\tz<1.$

\end{lem}
\begin{lem}\label{l2.4.}\hspace{-0.1cm}{\rm\bf 2.4(\cite{aj}).}\quad
Let $\wz\ge 0$ be locally integral on $\rz$, $1<r_i, s_i<\fz$, $1\le
p_i, q_i<\fz$ and suppose $T$ is a sublinear operator defined on
${\cal S}$ satisfying
$$\l(\int_\rz|Tf(x)|_{s_i}^{q_i}\wz(x)dx\r)^{1/q_i}\le
M_i\l(\int_\rz|f(x)|_{r_i}^{p_i}\wz(x)dx\r)^{1/p_i}$$ for $i=0,1$ and
$f\in {\cal S}$. Then $T$ extends uniquely to a sublinear operator
on $L_\wz^p(l^r)$ and there is a constant $M_\tz$ such that
$$\l(\int_\rz|Tf(x)|_r^q\wz(x)dx\r)^{1/q}\le M_0^{1-\tz}M_1^\tz\l(\int_\rz|f(x)|_r^{p}\wz(x)dx\r)^{1/p}$$
where $(1/p,1/q,1/s,1/r)=(1-\tz)(1/p_0,1/q_0, 1/s_0,
1/r_0)+\tz(1/p_1,1/q_1, 1/s_1, 1/r_1),\quad 0<\tz<1.$
\end{lem}
We define the dyadic maximal operator $M^\triangle_{V,\tz} f(x)$ as
follows
$$M_{V,\tz}^\triangle f(x):=\dsup_{x\in   Q(dyadic\ cube)}\dfrac
1{\psi_\tz(Q)|Q|}\dint_Q|f(x)|\,dx,$$   where  $\psi_\tz(Q)=(1+r/\max_Q\rho(x))^\tz$, $r$ is side-length of $Q$ and
 $\tz> 0$.
\begin{lem}\label{l2.5.}\hspace{-0.1cm}{\rm\bf 2.5.}\quad
Let $f$ be a  locally integrable function on $\rz$, $\lz>0$, and
$\Omega_\lz=\{x\in \rz:\ M^\triangle_{V,\tz}f(x)>\lz\}$. Then
$\Omega_\lz$ may be written as a disjoint union of dyadic cubes
$\{Q_j\}$ with
\begin{enumerate}
\item[(i)] $\lz<(\psi_\tz(Q_j)|Q_j|)^{-1}\dint_{Q_j}|f(x)|\,dx,$
\item[(ii)] $(\psi_\tz(Q_j)|Q_j|)^{-1}\dint_{Q_j}|f(x)|\,dx\le
(4n)^\tz 2^n\lz,$ for each cube $Q_j$. This has the immediate consequences:
\item[(iii)] $|f(x)|\le\lz \ {\rm for}\ a.e\ x\in
\rz\setminus\bigcup_jQ_j$
 \item[(iv)]$|\Omega_\lz|\le
\lz^{-1}\dint_\rz|f(x)|\,dx.$
\end{enumerate}
\end{lem}
The proof follows from the same argument of Lemma 1 in page 150 of
\cite{s}.

\begin{thm}\label{t2.1.}\hspace{-0.1cm}{\rm\bf 2.1.}\quad
Let $1<r<\fz$ and $\tz>0$.
\begin{enumerate}
\item[(a)] If $1\le p<\fz$, $\wz\in A_p^{\rho,\tz}$, $\eta=p_0\tz_0$ where $p_0=4(l_0+1)^5(p+(\frac {r+1}2)')$ and $\tz_0=p((3\tz+n)p+
(l_0+1)n)$, there is a constant $C_{r,p,\tz,l_0,C_0}$ such that
$$\wz(\{x\in\rz:\ |M_{V,\eta}f(x)|_r>\az\}|\le
C\az^{-p}\dint_\rz|f(x)|_r^p\wz(x)dx.\eqno(2.1)$$
\item[(b)] If $1<p<\fz$, $\wz\in A_p^{\rho,\tz}$ and $\eta$ be same as  above, there is a
constant $C_{r,p,\tz,l_0,C_0}$ such that
$$\dint_\rz|M_{V,\eta}f(x)|_r^p\wz(x)dx\le
C\az^{-p}\dint_\rz|f(x)|_r^p\wz(x)dx.\eqno(2.2)$$
\end{enumerate}
\end{thm}
{\it Proof.} \quad Observe first that (2.2) for the case $r=p$ is easy
consequence of Lemma 2.2 since $\eta>r'\tz$,
$$\begin{array}{cl}
\dint_\rz
|M_{V,\eta}f(x)|_r^r\wz(x)dx&=\dsum_k\dint_\rz|M_{V,\eta}f_k(x)|^r\wz(x)dx\\
&\le C\dsum_k\dint_\rz|f_k(x)|^r\wz(x)dx\\
&= C\dsum_k\dint_\rz|f_k(x)|_r^r\wz(x)dx.
\end{array}\eqno(2.3)$$
Now suppose $r>p$, $\wz\in A_p^{\rho,\tz}$ and $\az>0$. As usual, we
can assume  that $f\in C_0^\fz$. Let $\tz_1=\tz(l_0+1)$. From Lemma 2.5, we yields a
sequence of non-overlapping cube $\{Q_j\}$ such that
$$|f(x)|_r\le\az,\quad x\not\in\oz=\bigcup_{j=1}^\fz Q_j,\eqno(2.4)$$
$$\az<\dfrac1{\psi_{\tz_1}(Q_j)|Q_j|}\dint_{Q_j}|f(x)|_rdx\le
2^n(4n)^{\tz_1}\az,\quad j=1,2,\cdots.\eqno(2.5)$$ Let $f=f'+f''$ where
$f'=\{f_k'\}, f'_k(x)=f_k(x)\chi_{\rz\setminus\oz}(x)$. Then
$$|M_{V,\eta}f(x)|_r\le
|M_{V,\eta}f'(x)|_r+|M_{V,\eta}f''(x)|_r.$$ From this, (2.1) will
follow if we show that
$$\wz(\{x\in\rz:\ |M_{V,\eta}f'(x)|_r>\az\}|\le
C\az^{-p}\dint_\rz|f(x)|_r^p\wz(x)dx.\eqno(2.6)$$ and
$$\wz(\{x\in\rz:\ |M_{V,\eta}f''(x)|_r>\az\}|\le
C\az^{-p}\dint_\rz|f(x)|_r^p\wz(x)dx.\eqno(2.7)$$ Since $\wz\in
A_r^{\rho,\tz}$ by (i) of Proposition 2.1, from (2.3) and (2.4), we
then have
$$\wz(\{x\in\rz:\ |M_{V,\eta}f'(x)|_r>\az\}|\le
C\az^{-r}\dint_\rz|f(x)|_r^r\wz(x)dx\le
C\az^{-p}\dint_\rz|f(x)|_r^p\wz(x)dx.$$ Thus, (2.6) is proved. To
prove (2.7), define $\bar f=\{\bar f_k\}$ by
$$\bar f_k(x)=\dfrac1{\psi_{\tz_1}(Q_j)|Q_j|}\dint_{Q_j}|f_k(y)|dy, \quad
x\in Q_j,\ j=1,2,\cdots,$$ zero,  otherwise. Let $\bar Q_j=2nQ_j$.
We now claim that  for any $x\in\bar\oz=\bigcup_j\bar Q_j$,
$$M_{V,\eta}f_k''(x)\le CM_{V,\bar\eta}\bar f_k(x),\quad \forall
k,$$ where $\bar\eta=\eta/2(l_0+1)^2$.

In fact, $\forall\ x\not\in \bar \oz$, and cube $Q\ni x$, if
$Q_j\bigcap Q\not=\O$, then $Q_j\subset \bar Q=4nQ$, hence
$$\begin{array}{cl}
\dfrac1{\Psi_{\eta}(Q)|Q|}\dint_Q|f''_k(x)|dx&=
\dfrac1{\Psi_{\eta}(Q)|Q|}\dsum_{j}\dint_{Q_j\bigcap
Q}|f_k(x)|dx\\
&\le \dfrac1{\Psi_{\eta}(Q)|Q|}\dsum_{Q_j\subset\bar Q}\dint_{Q_j}|f_k(x)|dx\\
&\le \dfrac1{\Psi_{\eta}(Q)|Q|}\dsum_{Q_j\subset\bar Q}\psi_{\tz_1}(Q_j)\dint_{Q_j}\bar f_k(x)dx\\
&\le C\dfrac{\Psi_{\tz_2}(\bar Q)}{\Psi_{\eta}(Q)|Q|}\dint_{\bar Q}\bar f_k(x)dx\\
&\le CM_{V,\bar\eta}\bar f_k(x),
\end{array}$$
where $\tz_2=\tz_1(l_0+1)=\tz(l_0+1)^2$.

By the claim above, it is easy to see that (3.8) will follow if we
show
$$\wz(\bar\oz)\le C\az^{-p}\dint_\rz|f(x)|_r^p\wz(x)dx.\eqno(2.8)$$
and
$$\wz(\{x\in\rz:\ |M_{V,\bar\eta}\bar f(x)|_r>\az\}|\le
C\az^{-p}\dint_\rz|f(x)|_r^p\wz(x)dx.\eqno(2.9)$$ If $p>1$, by
(2.5), we then have
$$\begin{array}{cl}
\wz(\bar Q_j)&=\dint_{\bar Q_j}\wz(x)dx\le
\dfrac{\az^{-p}}{(\psi_{\tz_1}(Q)|Q|)^p}
\l(\dint_{Q_j}|f(x)|_r\r)^p\dint_{\bar Q_j}\wz(x)dx\\
&\le \az^{-p}
\l(\dint_{Q_j}|f(x)|^p_r\wz(x)dx\r)\l(\dfrac1{(\Psi_{\tz}(Q)|Q|)}\dint_{
Q_j}\wz^{-1/(p-1)}(x)dx\r)^{p-1}\\
&\qquad\qquad\times
\l(\dfrac1{(\Psi_{\tz}(Q)|Q|)}\dint_{\bar Q_j}\wz(x)dx\r)\\
&\le \az^{-p} \dint_{Q_j}|f(x)|^p_r\wz(x)dx,
\end{array}\eqno(2.10)$$
since $\wz\in A_p^{\rho,\tz}$.

 A similar
argument shows that (2.10) holds also if $p=1$. Hence, (2.8) follows
from (2.10) upon summing over $j$. Note that $|\bar f(x)|_r\le
2^n(4n)^{\tz_1}\az$, and since $|\bar f(x)|_r$ is supported in $\oz$, using Lemma
2.2, we obtain
$$\wz(\{x\in\rz:\ |M_{V,\bar\eta}\bar f(x)|_r>\az\}|\le
C\az^{-r}\dint_\rz|\bar f(x)|_r^r\wz(x)dx\le C\dint_\oz\wz(x)dx$$
which together with (2.10) yields (2.9) as required. This complete
the proof (2.1) in the case $r\ge p$. If $r>p>1$, by (iii) of
Proposition 2.1, we know that for $\wz\in A_p^{\rho,\tz}$, there
exist constants $p_1, p_2,\tz_3$(depending only on $\wz$) $(r+1)/2<p_1<p<p_2<r$ and $ \tz_3\le\tz_0$
so that (2.1) holds with $\wz\in A_{p_1}^{\tz_3}$ and $\wz\in
A_{p_2}^{\tz}$ respectively. Obviously, $\bar\eta>2p_1'\tz_3$, Lemmas
2.2 and 2.3 yields (2.2) for $r>p>1$.

Suppose now that $p>r$ and $\wz\in A_p^{\rho,\tz}$. By (iii) of
Proposition 2.1, there exist  constants $\tz_4\le\tz_0$ and $1<r_0<p$ such that $\wz\in A_q^{\rho,\tz_4}, q\ge
p/r_0$. In particular, (i) of Proposition 2.1 yield $\wz(x)>0$ a.e.
and $\wz(x)^{1-q'}\in A_{q'}^{\rho,\tz_4}$ so that by Lemma 2.2, for
any nonnegative function $\|\vz\|_{L^{q'}_\wz}\le 1$, we then have
$$\dint_\rz|M_{V,\eta_1}(\vz\wz)(x)|^{q'}\wz(x)^{1-q'}dx\le
C_q\dint_\rz|\vz(x)|^{q'}\wz(x)dx=C_q,$$ where
$\eta_1=\bar\eta/(l_0+1)^3>q\tz_4$ and hence
$$\begin{array}{cl}
\dint_\rz|M_{V,\bar\eta}f(x)|_r^r\vz(x)\wz(x)dx&\le
C\dint_\rz|f(x)|_r^r[M_{V,\eta_1}(\vz\wz)(x)/\wz^{1/q}(x)]\wz^{1/q}(x)dx\\
&\le C\l(\dint_\rz|f(x)|_r^{rq}\wz(x)dx\r)^{1/q}.
\end{array}\eqno(2.11)$$
In the first inequality of (2.11), we used the following fact that
for any nonnegative measurable functions $f,g$, and $q>1$, we have
$$\dint_\rz (M_{V,\bar\eta}f)^qgdx\le C\dint_\rz f^q
(M_{V,\eta_1}g)dx.\eqno(2.12)$$ Taking the supremum in (2.11) over
such $\vz$ then yields (2.2) for $1<r\le r_0$ upon taking $q=p/r$,
and this together with the case $p=r$ provided in (2.3) yields (3.3)
for $r_0<r<p$ by application of Lemma 2.4. Thus, the proof of (a)
and (b) is complete.

It remains to prove (2.12), let $\eta_2=\eta_1(l_0+1)=\bar\eta/(l_0+1)^2$, we shall begin by proving
$$\dint_\rz (M^\triangle_{V,\eta_2}f)^qgdx\le C\dint_\rz f^q
(M_{V,\eta_1}g)dx.\eqno(2.13)$$ We do this follows: Hold $g$ fixed,
and look at the mapping $T: f\to M^\triangle_{V,\eta_2}f$. Then
(2.13) says that $T$ is bounded from $L^q(\rz,M_{V,\eta_1}g(x)dx)$
to $L^q(\rz,g(x)dx)$. Clearly, $T$ is bounded from
$L^\fz(\rz,M_{V,\eta_1}g(x)dx)$ to $L^\fz(\rz,g(x)dx)$. If we can
show that $T$ is weak (1,1) type, then (2.13) holds by the
Marcinkiewicz interpolation theorem.

Lemma 2.1 shows that $\{x\in\rz:\
M^\triangle_{V,\eta_2}f(x)>\lz\}=\bigcup_j Q_j$, where the $Q_j$ are
pairwise disjoint cubes satisfying the condition
$$\lz\le \dfrac1{\psi_{\eta_2}(Q_j)|Q_j|}\dint_{Q_j}f(x)dx\le
2^n(4n)^{\eta_2}\lz.$$ Then
$$\begin{array}{cl}
\dint_{Q_j}g(y)dy&\le\dint_{Q_j}g(y)dy\dfrac{\lz^{-1}}{\psi_{\eta_2}(Q_j)|Q_j|}\dint_{Q_j}f(x)dx\\
&\le C\lz^{-1}\dint_{Q_j}f(x)\l[\dfrac1{\Psi_{\eta_1}(Q_j)|Q_j|}\dint_{Q_j}g(y)dy\r]dx\\
&\le C\lz^{-1}\dint_{Q_j}f(x)M_{V,\eta_1}g(x)dx.
\end{array}$$
Summing over $j$, we obtain
$$\dint_{\{x\in\rz:\
(M^\triangle_{V,\eta_2}f)(x)>\lz\}}g(y)dy\le C\dint_\rz
f(x)M_{V,\eta_1}g(x)dx,$$ Thus, (2.13) holds. To complete the
proof (2.12), we first define
$$M'_{V,\eta_3}f(x)=\dsup_{r>0}\dfrac1{(1+r/\rho(x))^{\eta_3}|Q|}\dint_{Q(x,r)}|f(y)|dy.$$
Obviously, $(4n)^{\bar\eta}C_0M'_{V,\eta_3}f(x)\ge M_{V,\bar\eta}f(x)$, where $\eta_3=\bar\eta/(l_0+1)=\eta_2(l_0+1)$.

Hence, to end the proof, it will suffice to show that
$$\{x\in\rz:\
M'_{V,\eta_3}f(x)>c_0\lz\}\subset\bigcup_j 2Q_j,\eqno(2.14)$$ where
$c_0=C_0^2 4^{l_0+1+n}(4n)^{\bar\eta}$.

Fix $x\in \bigcup_j 2Q_j$ and let $Q$ be any cube centered at $x$.
Let $r$ denote the side length of $Q$, and choose $k\in\zz$ such
that $2^{k-1}\le r<2^k$. Then $Q$ intersects $m(\le 2^n)$ dyadic cubes
with sidelengh $2^k$; call them $R_1=R_1(x_1, 2^k),R_2=R_2(x_2,
2^k),\cdots, R_m=R_m(x_m, 2^k)$. Non of these cubes is contained in
any of the $Q_j's$, for otherwise we would have $x\in
\bigcup_j(2Q_j)$. Hence
$$\begin{array}{cl}
\dfrac1{(1+r/\rho(x))^{\eta_3}|Q|}\dint_{Q(x,r)}|f(y)|dy&=
\dfrac1{(1+r/\rho(x))^{\eta_3}|Q|}\dsum_{i=1}^m\dint_{Q\bigcap
R_i}|f(y)|dy\\
&\le
\dsum_{i=1}^m\dfrac{C_04^{l_0+1}2^{kn}}{(1+2^k/\max_Q\rho(x))^{\eta_2}|Q||R_i|}\dint_{
R_i}|f(y)|dy\\
&\le 2^n 4^{l_0+1}C_0 m\lz\le 4^{l_0+1+n}C_0\lz.
\end{array}$$
Thus, (2.14) holds, so (2.12) is proved.\hfill$\Box$

\begin{center} {\bf 3. Extrapolation theorems}\end{center}
In this section, ${\cal F}$ will denote a family of order pairs of
non-negative, measurable function $(f,g)$. If we say that for $p$,
$0<p<\fz$, and $\wz\in A_\fz^{\rho,\fz}=\bigcup_{p\ge
1}A_p^{\rho,\fz}$.
$$\dint_\rz f(x)^p\wz(x)dx\le C\dint_\rz g(x)^p\wz(x),\quad (f,g)\in
{\cal F},$$ we mean that this inequality holds for any $(f.g)\in
{\cal F}$ such that the left-hand side is finite, and that the
constant $C$ depends only upon $p$ and the $ A_\fz^{\rho,\fz}$
constant of $\wz$. We will make similar abbreviated statements
involving Lorentz spaces. For vector-valued inequalities we will
consider sequences $\{(f_j,g_j)\}$, where each pair $(f_j,g_j)$ is
contained in ${\cal F}$.

In addition, we will use following classes: given a pair of
operators $T$ and $S$, let ${\cal F}(T,S)$ denote the family of
pairs of functions $(|Tf|, |Sf|)$, where $f$ lies in the common
domain of $T$ and $S$, and the left-hand side of the corresponding
inequality is finite.  To achieve this, the function $f$ may be
restricted in some other way, e.g. $f\in C_0^\fz$. In this case we
may indicate this by writing ${\cal F}(|Tf|,|Sf|: f\in C_0^\fz)$.

We can now state our main results in this paper.
\begin{thm}\label{t3.1.}\hspace{-0.1cm}{\rm\bf 3.1.}\quad
Given a family ${\cal F}$, suppose that for some $p_0,\ 0<p_0<\fz$,
and for every weight $\wz\in A_{\fz}^{\rho,\fz}$,
$$\dint_\rz f(x)^{p_0}\wz(x)dx\le C\dint_\rz g(x)^{p_0}\wz(x),\quad (f,g)\in
{\cal F}.\eqno(3.1)$$Then:

For all $0<p<\fz$ and $\wz\in A_{\fz}^{\rho,\fz}$
$$\dint_\rz f(x)^{p}\wz(x)dx\le C\dint_\rz g(x)^{p}\wz(x)dx,\quad (f,g)\in
{\cal F}.\eqno(3.2)$$

For all $0<p<\fz,\ 0<s\le\fz$ and $\wz\in A_{\fz}^{\rho,\fz}$
$$\|f\|_{L^{p,s}(\wz)}\le C\|g\|_{L^{p,s}(\wz)},\quad (f,g)\in
{\cal F}.\eqno(3.3)$$

For all $0<p,q<\fz$ and $\wz\in A_{\fz}^{\rho,\fz}$
$$\l\|\l(\dsum_j(f_j)^q\r)^{\frac 1q}\r\|_{L^{p}(\wz)}\le C\l\|\l(\dsum_j(g_j)^q\r)^{\frac 1q}\r\|_{L^{p}(\wz)},\quad
\{(f_j,g_j)\}_j\subset {\cal F}.\eqno(3.4)$$

For all $0<p,q<\fz$, $0<s\le\fz$, and $\wz\in A_{\fz}^{\rho,\fz}$
$$\l\|\l(\dsum_j(f_j)^q\r)^{\frac 1q}\r\|_{L^{p,s}(\wz)}\le C\l\|\l(\dsum_j(g_j)^q\r)^{\frac 1q}\r\|_{L^{p,s}(\wz)},\quad
\{(f_j,g_j)\}_j\subset {\cal F}.\eqno(3.5)$$
\end{thm}
Our second main result shows that we can also extrapolate from an
initial Lorentz space inequality.
\begin{thm}\label{t3.2.}\hspace{-0.1cm}{\rm\bf 3.2.}\quad
Given a family ${\cal F}$, suppose that for some $p_0, 0<p_0<\fz$,
and for every weight $\wz\in A_\fz^{\rho,\fz}$,
$$\|f\|_{L^{p_0,\fz}(\wz)}\le C\|g\|_{L^{p_0,\fz}(\wz)},\quad (f,g)\in
{\cal F}.\eqno(3.6)$$ For all $0<p<\fz$ and $\wz\in
A_{\fz}^{\rho,\fz}$
$$\|f\|_{L^{p,\fz}(\wz)}\le C\|g\|_{L^{p,\fz}(\wz)},\quad (f,g)\in
{\cal F}.\eqno(3.7)$$
\end{thm}

Our third main result is a generalization of the $A_p$ extrapolation
theorem of Rubio de Francia.

\begin{thm}\label{t3.3.}\hspace{-0.1cm}{\rm\bf 3.3.}\quad
Fix $\gz\ge 1$ and $r$, $\gz<r<\fz$. If $T$ is a bounded operator on
$L^r(\wz)$ for any $\wz\in A_{r/\gz}^{\rho,\fz}$, with operator norm
depending only the $A_{r/\gz}$ constant of $\wz$, then $T$ is
bounded on $L^p(\wz)$, $\gz<p<\fz$, for any $\wz\in
A_{p/\gz}^{\rho,\fz}$.
\end{thm}
As a consequence of Theorem 3.3, we have the following result.
\begin{cor}\label{c3.1.}\hspace{-0.1cm}{\rm\bf 3.1.}\quad
Fix $\gz\ge 1$. Let $\gz<p, q<\fz$ and $T$ satisfy the conditions in
Theorem 3.3. Then for any $\wz\in A_{p/\gz}^{\rho,\fz}$ such that
$$\l\|\l(\dsum_j|Tf_j|^q\r)^{\frac 1q}\r\|_{L^{p}(\wz)}\le C\l\|\l(\dsum_j|f_j|^q\r)^{\frac 1q}\r\|_{L^{p}(\wz)}.$$
\end{cor}

We shall adapt a similar argument in \cite{cmp} for proving Theorems
3.1 and 3.2, and prove Theorem 3.3 by using an argument in \cite{d}.
We first give the proof of Theorem 3.1.

\subsection*{ 3.1.\ Proof of inequality (3.2)}
We prove this inequality in two steps.

{\it Step 1}:\ We first show that hypothesis (3.1) is equivalent to
the family of weighted inequalities with $A_1^{\rho,\fz}$ weights.
\begin{prop}\label{p3.1.}\hspace{-0.1cm}{\rm\bf 3.1.}\quad
Hypothesis (3.1) of Theorem 3.1 is equivalent to the following: for
all $0<q<p_0,\ \wz\in A_1^{\rho,\fz}$, and $(f,g)\in{\cal F}$,
$$\dint_\rz f(x)^{q}\wz(x)dx\le C\dint_\rz g(x)^{q}\wz(x)dx.\eqno(3.8)$$
\end{prop}
{\bf Proof of Proposition 3.1.}\ We will prove that (3.1) implies
(3.8). If (3.2) is proved, then the converse is proved. Fix
$(f,g)\in{\cal F}$. Without loss of generality, we can assume that
$g\in L^q(\wz)$ and $\|f\|_{L^q(\wz)}>0$. Let $s=p_0/q$. Since
$\wz\in A_1^{\rho,\fz}$, so there is a $\tz>0$ such that $\wz\in
A_1^{\rho,\tz}\subset A_{s'}^{\rho,\tz}$, and $M_{V,s\tz}$ is
bounded on $L^{s'}(\wz)$ by Lemma 2.2, that is,
$$\|M_{V,s\tz}h\|_{L^{s'}(\wz)}\le A\|h\|_{L^{s'}(\wz)},$$ for some $A>0$.
For $h\in L^{s'}(\wz),\ h\ge 0$, we apply the algorithm of Rubio de
Francia to define
$${\cal R}h(x)=\dsum_{k=0}^\fz\dfrac{M^k_{V,s\tz}h(x)}{(2A)^k},$$
where $M^k_{V,s\tz}$ is the operator $M_{V,s\tz}$ iterated $k$ times
if $k\ge 1$, and for $k=0$ is just the identity. From the definition
of ${\cal R}$, it easy to see that:
\begin{enumerate}
\item[(a)] $h(x)\le {\cal R}h(x)$.
\item[(b)]$\|{\cal R}h\|_{ L^{s'}(\wz)}\le 2\|h\|_{ L^{s'}(\wz)}$.
\item[(c)]$M_{V,s\tz}({\cal R}h)(x)\le 2A {\cal R}h(x)$, so ${\cal
R}h(x)\in A_1^{\rho,s\tz}$ with constant independent of $h$.
\end{enumerate}
Since $f, g\in L^{s'}(\wz)$ and have positive norms, from (b), we
then have
$$H(x)={\cal R}\l(\l(\dfrac f{\|f\|_{ L^{s'}(\wz)}}\r)^{\frac{q}{s'}}
\l(\dfrac g{\|g\|_{ L^{s'}(\wz)}}\r)^{\frac{q}{s'}}\r)(x)\in
L^{s'}(\wz).$$ By (a),
$$\l(\dfrac f{\|f\|_{ L^{s'}(\wz)}}\r)^{\frac{q}{s'}}\le H(x),\
\l(\dfrac g{\|g\|_{ L^{s'}(\wz)}}\r)^{\frac{q}{s'}}\le
H(x),\eqno(3.9)$$ So $H(x)>0$ whenever $f(x)>0$. Further, $H$ is
finite a.e. on the set where $\wz>0$ because $h\in L^{s'}(\wz)$.
Hence,
$$\dint_\rz f(x)^q\wz(x)dx\le \l(\dint_\rz
f(x)^{p_0}H(x)^{-s}\wz(x)dx\r)^{\frac 1s} \l(\dint_\rz
H(x)^{s'}\wz(x)dx\r)^{\frac 1{s'}}:=I\cdot II.$$ Obviously, $II\le 4$
by (b).

To estimate $I$, since $\wz\in A_1^{\rho,\tz}\subset
A_1^{\rho,s\tz}$, and $H\in A_1^{\rho,s\tz}$ by (c), so
$wH^{-s}=wH^{1-(1+s)}\in A_{1+s}^{\rho,s\tz}\subset
A_{\fz}^{\rho,\fz} $ by (v) of Proposition 2.1.  on the other hand,
by (3.9), we have
$$\dint_\rz
f(x)^{p_0}H(x)^{-s}\wz(x)dx\le
\|f\|^{\frac{qs}{s'}}_{L^s(\wz)}\dint_\rz
f(x)^{p_0-\frac{qs}{s'}}\wz(x)dx=\|f\|^{qs}_{L^s(\wz)}<\fz.$$ So,
we can use (3.1); by (3.9), we get
$$I\le\l(\dint_\rz
g(x)^{p_0}H(x)^{-s}\wz(x)dx\r)^{\frac 1s}\le C\dint_\rz
g(x)^p\wz(x)dx.$$ By I and II, we obtain the desired result.

\medskip

{\it Step 2}:\  We now show that for all $0<p<\fz$ and for every
$\wz\in A_\fz^{\rho,\fz}$, (3.2) holds. Fix $0<p<\fz$ and $\wz\in
A_\fz^{\rho,\fz}$. Assume that $(f,g)\in{\cal F}$ with  $f\in
L^p(\wz)$ and $g\in L^p(\wz)$. By (i) of Proposition 2.1, we know that
$A_{p_1}^{\rho,\tz}\subset A_{p_2}^{\rho,\tz}$ if $1\le p_1\le p_2$,
there exist $\tz>0$ and $0<q<\min\{p,p_0\}$ such that $\wz\in
A_{p/q}^{\rho,\tz}$. Let $r=p/q>1$. Since $\wz\in A_{r}^{\rho,\tz}$,
then $\wz^{1-r'}\in A_{r'}^{\rho,\tz}$ by (ii) of Proposition 2.1.
Given $h\in L^{r'}(\wz^{1-r'}),\ h\ge 0$, we use the algorithm of
Rubio de Francia to define
$${\cal R}h(x)=\dsum_{k=0}^\fz\dfrac{M^k_{V,r\tz}h(x)}{(2B)^k},$$
where $B$ is the operator norm of $M_{V,r\tz}$ on
$L^{r'}(\wz^{1-r'})$; this is finite since $\wz^{1-r'}\in
A_{r'}^{\rho,\tz}$. Then
\begin{enumerate}
\item[(a)] $h(x)\le {\cal R}h(x)$.
\item[(b)]$\|{\cal R}h\|_{ L^{r'}(\wz^{1-r'})}\le 2\|h\|_{  L^{r'}(\wz^{1-r'})}$.
\item[(c)]$M_{V,sr}({\cal R}h)(x)\le 2B {\cal R}h(x)$, so ${\cal
R}h(x)\in A_1^{\rho,r\tz}$ with constant independent of $h$.
\end{enumerate}

By duality
$$\|f\|^q_{L^p(\wz)}=\|f^q\|_{L^r(\wz)}=\dsup_{\|h\|_{L^{r'}(\wz)}\le 1}\dint_\rz
f(x)^qh(x)\wz(x).$$ Fix such a function $h\ge 0$. Then $h\wz\in
L^{r'}(\wz^{1-r'})$ and
$\|h\wz\|_{L^{r'}(\wz^{1-r'})}=\|h\|_{L^{r'}(\wz)}=1$. By (c),
${\cal R}(h\wz)\in A_1^{\rho,r\tz}$. By (a) and Proposition 3.1, we
then have
$$\dint_\rz f(x)^qh(x)\wz(x)dx\le \dint_\rz f(x)^q{\cal
R}(h\wz)(x)dx\le C\dint_\rz g(x)^q{\cal R}(h\wz)(x)dx,$$ provided
that the middle term is finite, this is obvious.

The same argument also holds for $g$ instead of $f$. Hence,
$$\dint_\rz f(x)^qh(x)\wz(x)dx\le C\dint_\rz g(x)^q{\cal R}(h\wz)(x)dx\le
C\|g\|_{L^p(\wz)}^q.$$ From this, we obtain the desired result.\hfill$\Box$

\subsection*{ 3.2.\ Proof of inequality (3.3)}
We  need two lemmas. We first give a result about the operator
$M_{\wz}$ defined by
$$M_{\wz}(f)(x)=\dsup_{x\in B}\dfrac
1{\wz(5B)}\dint_B|f(x)|\wz(x)dx.$$
\begin{lem}\label{l3.1.}\hspace{-0.1cm}{\rm\bf 3.1(\cite{t1}).}\quad
Let  $1\le p<\fz$ . If $\wz\in
A_\fz^{\rho,\fz}$ , then
$$\wz(\{x\in\rz:\ M_{\wz}f(x)>\lz\})\le C\l(\dfrac
{\|f\|_{L^p(\wz)}}\lz\r)^{p},\ \forall\lz>0,\ \forall f\in
L^p(\wz).$$ In particular, for $1<p\le fz$,
$$\|M_{\wz}f\|_{L^p(\wz)}\le C\|f\|_{L^p(\wz)}.$$
\end{lem}
Given two weights $u$ and $v$, we say that
$u\in A_{1}(v)$ if for every $x$, $M_{v}u(x)\le Cu(x).$
\begin{lem}\label{l3.2.}\hspace{-0.1cm}{\rm\bf 3.2.}\quad
If $\wz_1\in A_{p}^{\rho,\tz}$, $1\le p\le\fz$, and $\wz_2\in
A_1(\wz_1)$, then $\wz_1\wz_2\in A_p^{\rho,\tz p}$.
\end{lem}
{\it Proof.}\quad If $\wz_2\in A_1(\wz_1)$, then for any ball $B$
$$\begin{array}{cl}
\dfrac1{(\Psi_\tz(B))^{p^2}|B|}\dint_B\wz_1(x)\wz_2(x)dx&=\dfrac{\wz_1(5B)}{(\Psi_\tz(B))^{p^2}|B|}\dfrac1
{\wz_1(5B)} \dint_B\wz_2(x)\wz_1(x)dx\\
&\le C \dfrac{\wz_1(5B)}{(\Psi_\tz(B))^{p^2}|B|}{\rm ess}\dinf_B \wz_2\\
&\le C \dfrac{\wz_1(B)}{(\Psi_\tz(B))^p|B|}{\rm ess}\dinf_B \wz_2,
\end{array}$$
in the last inequality, we used the following fact (see \cite{t1})
$$\wz_1(5B)\le C(\Psi_\tz(B))^p\wz_1(B).$$
On the other hand,
$$\l(\dfrac1{|B|}\dint_B(\wz_1(x)\wz_2(x))^{-\frac1{p-1}}dx\r)^{p-1}\le
\l(\dfrac1{|B|}\dint_B\wz_1(x)^{-\frac1{p-1}}dx\r)^{p-1}({\rm ess}\dinf_B
\wz_2)^{-1}.$$ From two inequalities above, we get the desired
result.\hfill$\Box$

\medskip

{\it Proof of (3.3).}\quad Fix $p,s$, $\wz\in A_\fz^{\rho,\fz}$ and
$(f,g)\in{\cal F}$ with $f,g\in L^{p,s}(\wz)$. Fix $0<q<\min\{p,s\}$
and set $r=p/q>1$, $\wt r=s/q>1$. (If $s=\fz$, take $0<q<p$ and $\wt
r=\fz$.) Then
$$\|f\|_{L^{p,s}(\wz)}^q=\|f^q\|_{L^{r,\wt r}(\wz)}=\dsup_h\dint_\rz
f(x)^qh(x)\wz(x)dx,$$ where the supremum is taken over all $h\in
L^{r',\wt r}(\wz)$ with $h\ge 0$ and $\|h\|_{L^{r',\wt r'}}=1$. Fix
such a function $h$. Using the algorithm of Rubio de Francia to
define
$${\cal R}_\wz h(x)=\dsum_{k=0}^\fz\dfrac{M^k_{\wz}h(x)}{(2A_\wz)^k},$$
where $A_\wz$ is the operator norm of $M_\wz$ on $L^{r'\wt r}(\wz)$
endowed with norm equivalent to $\|\cdot\|_{L^{r',\wt r}(\wz)}$.
Since $M_\wz$ is bounded on $L^p(\wz)$ by Lemma 3.1, and by Marcinkiewicz
interpolation in the scale of Lorentz space, it is bounded on
$L^{r',\wt r}(\wz)$. Then,
\begin{enumerate}
\item[(a)] $h(x)\le {\cal R}_\wz h(x)$.
\item[(b)]$\|{\cal R}_\wz h\|_{ L^{r',\wt r}(\wz^{1-r'})}\le C\|h\|_{  L^{r',\wt r}(\wz^{1-r'})}=C$.
\item[(c)]$M_{V,s\tz}({\cal R}h)(x)\le 2A_\wz {\cal R}h(x)$, so ${\cal
R}_\wz h(x)\in A_1(\wz)$ with constant independent of $h$.
\end{enumerate}
By Lemma 3.2, $\wz{\cal R}_\wz h\in A_\fz^{\rho,\fz}$. As above,
(3.2) holds with exponent $q$ and the $A_\fz^{\rho,\fz}$ weight
$\wz{\cal R}_\wz h$. Thus,
$$\begin{array}{cl}
\dint_\rz f(x)^qh(x)\wz(x)dx&\le \dint_\rz f(x)^q{\cal
R}_\wz h(x)\wz(x)dx\le C\dint_\rz g(x)^q{\cal R}_\wz
h(x)\wz(x)dx\\
&\le C\|g^q\|_{L^{r,\wt r}(\wz)}\|{\cal R}_\wz h\|_{L^{r',\wt
r'}(\wz)}\le C\|g\|^q_{L^{r,\wt r}(\wz)},
\end{array}$$
since
$$\dint_\rz f(x)^q{\cal
R}_\wz h(x)\wz(x)dx\le \|f^q\|_{L^{r,\wt r}(\wz)}\|{\cal R}_\wz
h\|_{L^{r',\wt r'}(\wz)}\le C\|f\|^q_{L^{r,\wt r}(\wz)}<\fz.$$ Thus,
the desired inequality is obtained.\hfill$\Box$

\subsection*{ 3.3.\ Proof of inequalities (3.4) and (3.5)}
Fix $0<q<\fz$. It suffices to prove the vector-valued inequalities
only for finite sums by the monotone convergence theorem. Fix $N\ge
1$ and define
$$f_q(x)=\l(\dsum_{j=1}^N f_j(x)^q\r)^{\frac 1q},\
g_q(x)=\l(\dsum_{j=1}^N g_j(x)^q\r)^{\frac 1q},$$ where
$\{(f_j,g_j)\}_{j=1}^N\subset{\cal F}$.  Now form a new family
${\cal F}_q$ consisting of the pairs $(f_q,g_q)$. Then, for every
$\wz\in A_\fz^{\rho,\fz}$ and $(f_q,g_q)\in {\cal F}_q$, by (3.2) we
get
$$\|f_q\|_{L^q(\wz)}^q=\dsum_{j=1}^N\dint_\rz f_j(x)^q\wz(x)dx\le
C\dsum_{j=1}^N\dint_\rz g_j(x)^q\wz(x)dx=C\|g_q\|^q_{L^q(\wz)},$$
which implies that the hypotheses of Theorem 3.1 are fulfilled by
${\cal F}_q$ with $p_0=q$. Hence, by (3.2) and (3.3), for all,
$0<p<\fz,\ 0<s\le\fz,\ \wz\in A_\fz^{\rho,\fz}$, and $(f_q,g_q)\in
{\cal F}_q,$ $\|f_q\|_{L^p(\wz)}\le C\|g_q\|_{L^p(\wz)}$ and
$\|f_q\|_{L^{p,s}(\wz)}\le C\|g_q\|_{L^{p,s}(\wz)}$.\hfill$\Box$

\subsection*{ 3.4.\ Proof of Theorem 3.2}
Similar to the proof of Theorem 3.1, and adapting the same argument
of Theorem 2.2 in \cite{cmp}, we omit the details here.

\subsection*{ 3.4.\ Proof of Theorem 3.3}
We first need the following lemma, which is different from Lemma 2.2.
\begin{lem}\label{l3.3.}\hspace{-0.1cm}{\rm\bf 3.3(\cite{t1}).}\quad
Let $1\le p <\fz$ and suppose that $\wz\in A_p^{\rho,\tz}$. If $p<p_1<\fz$,
then
$$\dint_\rz|M_{V,\tz}f(x)|^{p_1}\wz(x)dx\le
C\dint_\rz|f(x)|^{p_1}\wz(x)dx.$$
\end{lem}

{\it Proof}.\quad  We only consider the case $\gz=1$, another case $\gz>1$
is similar.  We first show that if $1<q<r$ and $\wz\in
A_1^{\rho,\fz}$ then $T$ is bounded on $L^q(\wz)$. Without loss of
generality, we assume $\wz\in A_1^{\rho,\eta}$ for some $\eta>0$. By
(vi) of Proposition 2.1 the function $M_{V,\eta}^{(r-q)/(r-1)}$ is
in $A_1^{\rho,\eta}$, and $\wz (M_{V,\eta}f)^{q-r}\in
A_r^{\rho,\eta}$ by (v) of Proposition 2.1. Hence,
$$\begin{array}{cl}
\dint_\rz|Tf|^q\wz&=\dint_\rz |Tf|^q(M_{V,\eta}f)^{-(q-r)q/r}(M_{V,\eta}f)^{(q-r)q/r}\wz\\
&\le \l(\dint_\rz|Tf|^r\wz (M_{V,\eta}f)^{q-r}\r)^{q/r}\l(\dint_\rz
(M_{V,\eta}f)^q\wz\r)^{(r-q)/r}\\
&\le \l(\dint_\rz|f|^r\wz (M_{V,\eta}f)^{q-r}\r)^{q/r}\l(\dint_\rz
|f|^q\wz\r)^{(r-q)/r}\\
&\le C\dint_\rz |f|^q\wz,
\end{array}$$
the second inequality holds by our hypothesis on $T$ and by Lemma
3.3 (since $\wz\in A_1^{\rho,\eta}$), and the third inequality holds
since $|f(x)|\le M_{V,\eta}f(x)$ a.e. for any $\eta\ge 0$, so
$M_{V,\eta}f(x)^{q-r}\le |f(x)|^{q-r}$ a.e.

Given any $1<p<\fz$ and $\wz\in A_p^{\rho,\tz}$, by (iii) of
Proposition 2.1 there exists  $q>1$ and $\tz_1\ge\tz$ such that
$\wz\in A_{p/q}^{\rho,\tz_1}$, hence we only need to prove that $T$
is bounded on $L^p(\wz)$ if $\wz\in
A_{p/q}^{\rho,\tz_1}$.

Fix $\wz\in A_{p/q}^{\rho,\tz_1}$. Then by duality there exists
$u\in L^{(p/q)'}(\wz)$ with norm $1$ such that
$$\l(\dint_\rz|Tf|^p\wz\r)^{q/p}=\dint_\rz|Tf|^q\wz u.$$
 For any $s>1$, $\wz u\le M_{V,\eta}((\wz u)^s)^{1/s}$ for any
 $\eta>0$ and $M_{V,\eta}((\wz u)^s)^{1/s}\in A_1^{\rho,\eta}$.
 Hence, by the first part of the proof,
$$\begin{array}{cl}
\dint_\rz|Tf|^q\wz u&\le \dint_\rz|Tf|^qM_{V,\eta}((\wz
u)^s)^{1/s}\\
&\le C\dint_\rz|f|^qM_{V,\eta}((\wz
u)^s)^{1/s}\\
&=C\dint_\rz|f|^q\wz^{q/p}M_{V,\eta}((\wz u)^s)^{1/s}\wz^{-q/p}\\
&\le C\l(\dint_\rz|f|^p\wz\r)^{q/p}\l(\dint_\rz
M_{V,\eta}((\wz u)^s)^{(p/q)'/s}\wz^{1-(p/q)'}\r)^{1/(p/q)'}\\
\end{array}$$
Since $\wz\in A_{p/q}^{\rho,\tz_1}$,  then $\wz^{1-(p/q)'}\in
A_{(p/q)'}^{\rho,\tz_1}$ by (ii) of Proposition 2.1. Therefore, if
take $s$ sufficient close to 1, then there exists $\tz_s$ such that
$\wz^{1-(p/q)'}\in A_{(p/q)'/s}^{\rho,\tz_s}$ by (iii) of
Proposition 2.1. If choosing $\eta=((p/q)'/s)'\tz_s$, then by Lemma
2.2 the second integral is dominated by
$$C\dint_\rz(\wz u)^{(p/q)'}\wz^{1-(p/q)'}=C.$$
The proof is complete.\hfill$\Box$

\begin{center} {\bf 4. Some applications}\end{center}
\subsection*{4.1. Schr\"odinger type operators}
Let $T$ be a Schr\"odinger type operators. From Theorem 3.1 in
\cite{t1} we know that for all $0<p<\fz$ and $\wz\in A^{\rho,\fz}_\fz$, for
any $\eta>0$, then there exists a constant $C$ depending only on
$\eta, p, q, C_0, l_0$ and the $A_\fz^{\rho,\fz}$ constant of $\wz$
such that
$$\|Tf\|_{L^p(\wz)}\le C\|M_{V,\eta}f\|_{L^p(\wz)}.$$
By applying Theorem 3.1 to the family ${\cal F}_\eta(|Tf|,
M_{V,\eta}f: f\in C_0^\fz)$, we obtain that

For all $0<p,q<\fz$ and $\wz\in A_{\fz}^{\rho,\fz}$
$$\l\|\l(\dsum_j|Tf_j|^q\r)^{\frac 1q}\r\|_{L^{p}(\wz)}\le C\l\|\l(\dsum_j(M_{V,\eta}f_j)^q\r)^{\frac 1q}\r\|_{L^{p}(\wz)},\quad
\{(f_j,g_j)\}_j\subset {\cal F}_\eta.\eqno(4.1)$$

For all $0<p,q<\fz$, $0<s\le\fz$, and $\wz\in A_{\fz}^{\rho,\fz}$
$$\l\|\l(\dsum_j|Tf_j|^q\r)^{\frac 1q}\r\|_{L^{p,s}(\wz)}\le C\l\|\l(\dsum_j(M_{V,\eta}f_j)^q\r)^{\frac 1q}\r\|_{L^{p,s}(\wz)},\quad
\{(f_j,g_j)\}_j\subset {\cal F}_\eta.\eqno(4.2)$$

If we combine them with Theorem 2.1, we have the following
inequalities:

If $1<q<\fz$, then for every $\wz\in A_1^{\rho,\fz}$, there exists a
constant $C$ depending only on $\eta, q, C_0, l_0$ and the
$A_1^{\rho,\fz}$ constant of $\wz$ such that
$$\l\|\l(\dsum_j|Tf_j|^q\r)^{\frac 1q}\r\|_{L^{1,\fz}(\wz)}
\le C\l\|\l(\dsum_j|f_j|^q\r)^{\frac 1q}\r\|_{L^1(\wz)},\eqno(4.3)$$

If $1<q<\fz$, and $1<p<\fz$, then for every $\wz\in A_p^{\rho,\fz}$,
there exists a constant $C$ depending only on $\eta, p, q, C_0, l_0$
and the $A_p^{\rho,\fz}$ constant of $\wz$ such that
$$\l\|\l(\dsum_j|Tf_j|^q\r)^{\frac 1q}\r\|_{L^{p}(\wz)}
\le C\l\|\l(\dsum_j(f_j)^q\r)^{\frac 1q}\r\|_{L^p(\wz)}.\eqno(4.4)$$
Let $T$ be a Schr\"odinger type operators as above. From Theorem 3.1
in \cite{t1} we have that for all $0<p<\fz$ and $\wz\in A_\fz$, for
any $\eta>0$, then there exists a constant $C$ depending only on
$\eta, p, q, C_0, l_0$ and the $A_\fz^{\rho,\fz}$ constant of $\wz$
such that
$$\|[b,T]f\|_{L^p(\wz)}\le C\|b\|_{BMO_\fz(\rho)}\|M_{V,\eta}(M_{V,\eta}f)\|_{L^p(\wz)}.$$
By applying Theorem 3.1 to the family ${\cal F}_\eta(|[b,T]f|,
M_{V,\eta}f: f\in C_0^\fz)$, we obtain that

For all $0<p,q<\fz$ and
$\wz\in A_{\fz}^{\rho,\fz}$
$$\l\|\l(\dsum_j|[b,T]f_j|^q\r)^{\frac 1q}\r\|_{L^{p}(\wz)}\le C\|b\|_{BMO_\fz(\rho)}\l\|\l(\dsum_j(M_{V,\eta}f_j)^q\r)^{\frac 1q}\r\|_{L^{p}(\wz)},\quad
\{(f_j,g_j)\}_j\subset {\cal F}_\eta.\eqno(4.5)$$

For all $0<p,q<\fz$, $0<s\le\fz$, and $\wz\in A_{\fz}^{\rho,\fz}$
$$\l\|\l(\dsum_j|[b,T]f_j|^q\r)^{\frac 1q}\r\|_{L^{p,s}(\wz)}\le C\|b\|_{BMO_\fz(\rho)}\l\|\l(\dsum_j(M_{V,\eta}f_j)^q\r)^{\frac 1q}\r\|_{L^{p,s}(\wz)},\
\{(f_j,g_j)\}_j\subset {\cal F}_\eta,\eqno(4.6)$$
where the new space $BMO_\tz(\rho)$ introduced in \cite {bhs1} as follows
$$\|f\|_{BMO_\tz(\rho)}=\dsup_{B\subset
\rz}\dfrac 1{\Psi_\tz(B)|B|}\dint_B|f(x)-f_B|dx<\fz,$$ where
$f_B=\frac 1{|B|}\int_Bf(y)dy$ and $\Psi_\tz(B)=(1+r/\rho(x_0))^\tz$, $B= B(x_0,r)$ and $\tz> 0$.
Let $BMO_\fz(\rho)$ denote $\bigcup_{\tz>0} BMO_\tz(\rho)$

 If we combine them
with Theorem 2.1, we have the following inequality: If $1<q<\fz$,
and $1<p<\fz$, then for every $\wz\in A_p^{\rho,\fz}$, there exists
a constant $C$ depending only on $\eta, p, q, C_0, l_0$ and the
$A_p^{\rho,\fz}$ constant of $\wz$ such that
$$\l\|\l(\dsum_j|[b,T]f_j|^q\r)^{\frac 1q}\r\|_{L^{p}(\wz)}
\le C\|b\|_{BMO_\fz(\rho)}\l\|\l(\dsum_j|f_j|^q\r)^{\frac
1q}\r\|_{L^p(\wz)}.\eqno(4.7)$$
We remark  that these inequalities (4.1)-(4.7) are all new.

Next we consider another class $V\in B_q$ for $n/2\le q$ for Riesz transforms associated to Schr\"odinger operators. Let $T_1=(-\triangle+V)^{-1}V,\
T_2= (-\triangle+V)^{-1/2} V^{1/2}$ and $T_3=(-\triangle+V)^{-1/2} \nabla$. By using Theorem 3.3 in \cite{t2} and  Corollary 3.3, we have
\begin{thm}\label{t4.1.}\hspace{-0.1cm}{\rm\bf 4.1.}\quad
Suppose $V\in B_q$ and $q\ge n/2$.  Then
\begin{enumerate}
\item[(i)] If $q'< p,r<\fz$ and $\wz\in A_{p/q'}^{\rho,\fz}$,
$$\| |T_1f|_r\|_{L^p(\wz)}\le C\||f|_r\|_{L^p(\wz)};$$
\item[(ii)] If $(2q)'< p,r<\fz$ and $\wz\in A_{p/(2q)'}^{\rho,\fz}$,
$$\| |T_2f|_r\|_{L^p(\wz)}\le C\||f|_r\|_{L^p(\wz)};$$
\item[(iii)] If $p_0'< p,r<\fz$ and $\wz\in A_{p/p_0'}^{\rho,\fz}$, where $1/p_0=1/q-1/n$ and $n/2\le q<n$,
$$\|| T_3f|_r\|_{L^p(\wz)}\le C\||f|_r\|_{L^p(\wz)}.$$
\end{enumerate}
\end{thm}

Let $T_1^*=V(-\triangle+V)^{-1}, T^*_2=  V^{1/2}(-\triangle+V)^{-1/2}$ and $T^*_3=\nabla(-\triangle+V)^{-1/2} $. By duality we can easily get the following results.
\begin{cor}\label{c4.1.}\hspace{-0.1cm}{\rm\bf 4.1.}\quad
Suppose $V\in B_q$ and $q\ge n/2$.  Then
\begin{enumerate}
\item[(i)] If $1< p,r< q$ and $\wz^{-\frac 1{p-1}}\in A_{p'/q'}^{\rho,\fz}$,
$$\| |T^*_1f|_r\|_{L^p(\wz)}\le C\||f|_r\|_{L^p(\wz)};$$
\item[(ii)] If $1< p,r<2q$ and $\wz^{-\frac 1{p-1}}\in A_{p'/(2q)'}^{\rho,\fz}$,
$$\| |T^*_2f|_r\|_{L^p(\wz)}\le C\||f|_r\|_{L^p(\wz)};$$
\item[(iii)] If $1< p,r< p_0$ and $\wz^{-\frac 1{p-1}}\in A_{p'/p'_0}^{\rho,\fz}$, where $1/p_0=1/q-1/n$ and $n/2\le q<n$,
$$\| |T^*_3f|_r\|_{L^p(\wz)}\le C\||f|_r\|_{L^p(\wz)}.$$
\end{enumerate}
\end{cor}
Let $T_1,\ T_2$ and $T_3$ be above. By using Theorem 4.5 in \cite{t2} and  Corollary 3.3, we have
\begin{thm}\label{t4.2.}\hspace{-0.1cm}{\rm\bf 4.2.}\quad
Suppose $V\in B_q$ and $q\ge n/2$. Let $b\in BMO_\fz(\rho)$. Then
\begin{enumerate}
\item[(i)] If $q'< p,r<\fz$ and $\wz\in A_{p/q'}^{\rho,\fz}$,
$$\| |[b,T_1]f|_r\|_{L^p(\wz)}\le C\|b\|_{BMO_\fz(\rho)}\||f|_r\|_{L^p(\wz)};$$
\item[(ii)] If $(2q)'< p,r<\fz$ and $\wz\in A_{p/(2q)'}^{\rho,\fz}$,
$$\| |[b,T_2]f|_r\|_{L^p(\wz)}\le C\|b\|_{BMO_\fz(\rho)}\||f|_r\|_{L^p(\wz)};$$
\item[(iii)] If $p_0'< p,r<\fz$ and $\wz\in A_{p/p_0'}^{\rho,\fz}$, where $1/p_0=1/q-1/n$ and $n/2\le q<n$,
$$\| |[b,T_3]f|_r\|_{L^p(\wz)}\le C\|b\|_{BMO_\fz(\rho)}\||f|_r\|_{L^p(\wz)}.$$
\end{enumerate}
\end{thm}

Let $T_1^*,\ T^*_2$ and $T^*_3$ be above. By duality we can easily get the following results.
\begin{cor}\label{c4.2.}\hspace{-0.1cm}{\rm\bf 4.2.}\quad
Suppose $V\in B_q$ and $q\ge n/2$.  Let $b\in BMO_\fz(\rho)$. Then
\begin{enumerate}
\item[(i)] If $1< p,r< q$ and $\wz^{-\frac 1{p-1}}\in A_{p'/q'}^{\rho,\fz}$,
$$\||[b, T^*_1]f|_r\|_{L^p(\wz)}\le C\|b\|_{BMO_\fz(\rho)}\||f|_r\|_{L^p(\wz)};$$
\item[(ii)] If $1< p,r<2q$ and $\wz^{-\frac 1{p-1}}\in A_{p'/(2q)'}^{\rho,\fz}$,
$$\| |[b,T^*_2]f|_r\|_{L^p(\wz)}\le C\|b\|_{BMO_\fz(\rho)}\||f|_r\|_{L^p(\wz)};$$
\item[(iii)] If $1< p,r<p_0$ and $\wz^{-\frac 1{p-1}}\in A_{p'/p'_0}^{\rho,\fz}$, where $1/p_0=1/q-1/n$ and $n/2\le q<n$,
$$\| |[b,T^*_3]f|_r\|_{L^p(\wz)}\le C\|b\|_{BMO_\fz(\rho)}\||f|_r\|_{L^p(\wz)}.$$
\end{enumerate}
\end{cor}
Finally, we consider the Littlewood-Paley $g$ function related to Schr\"odinger operators is defined by
$$g(f)(x)=\l(\dint_0^\fz\l|\frac d{dt}e^{-tL}(f)(x)\r|^2tdt\r)^{1/2},$$ and the commutator $g_b$ of $g$ with $b\in BMO(\rho)$ is defined by
$$g_b(f)(x)=\l(\dint_0^\fz\l|\frac d{dt}e^{-tL}((b(x)-b(\cdot))f)(x)\r|^2tdt\r)^{1/2}.$$
The maximal operator of the diffusion semi-group is defined by
$$T^*f(x)=\dsup_{t>0}|e^{-tL}f(x)|=\dsup_{t>0}\l|\dint_\rz k_t(x,y)f(y)dy\r|,$$
 and it's commutator
$$T_b^*f(x)=\dsup_{t>0}\l|\dint_\rz k_t(x,y)(b(x)-b(y))f(y)dy\r|,$$
where $k_t$ is the kernel of the operator $e^{-tL},\ t>0$.

By Combining Theorems  1 and 2 in \cite{bhs2} and Theorems 1.1 and 3.1 in \cite{t2} and  Corollary 3.3 together, we have
\begin{thm}\label{t4.3.}\hspace{-0.1cm}{\rm\bf 4.3.}\quad Let $b\in BMO_\fz(\rho)$ and $T, T_b^*$, $g$ and $g_b$ be as above.
\begin{enumerate}
\item[(i)] If   $1< p,r<\fz$,
$\wz\in A_p^{\rho,\fz}$, then there exists a constant $C$ such that
$$\||g(f)|_r\|_{L^p(\wz)}+\||T^* f|_r\|_{L^p(\wz)}
\le C \||f|_r\|_{L^p(\wz)}.$$
\item[(ii)]If   $1< p,r<\fz$,
$\wz\in A_p^{\rho,\fz}$, then there exists a constant $C$ such that
$$\||g_b(f)|_r\|_{L^p(\wz)}+\||T^*_b f|_r\|_{L^p(\wz)}
\le C \|b\|_{BMO_\fz(\rho)}\||f|_r\|_{L^p(\wz)}.$$
\end{enumerate}
\end{thm}

\subsection*{4.2. Pseudo-differential operators}

Let $m$ be real number. Following \cite{t}, a symbol in
$S_{1,\dz}^m$ is a smooth function $\sz(x,\xi)$ defined on
$\rz\times\rz$ such that for all multi-indices $\az$ and $\bz$ the
following estimate holds:
$$|D_x^\az D_\xi^\bz\sz(x,\xi)|\le
C_{\az,\bz}(1+|\xi|)^{m-|\bz|+\dz|\az|},$$ where $ C_{\az,\bz}>0$
is independent of $x$ and $\xi$. A symbol in $S_{1,\dz}^{-\fz}$ is
one which satisfies the above estimates for each real number $m$.

The operator $T$ given by
$$Tf(x)=\dint_\rz\sz(x,\xi)e^{2\pi ix\cdot \xi}\hat f(\xi)\,d\xi$$
is called a pseudo-differential operator with symbol
$\sz(x,\xi)\in S_{1,\dz}^m$, where $f$ is a Schwartz function and
$\hat f$ denotes the Fourier transform of $f$. As usual,
$L^m_{1,\dz}$ will denote the class of pseudo-differential
operators with symbols in $S_{1,\dz}^m$.

We in \cite{t3} studied weighted inequalities for a class of pseudo-differential operators with symbols in $S_{1,\dz}^0$ with $0<\dz<1$.
More precisely, we  have the following result.
\begin{lem}\label{l4.1.}\hspace{-0.1cm}{\rm\bf 4.1.}\quad
Let $T$ be a pseudo-differential operators with symbols in $S_{1,\dz}^0$ with $0<\dz<1$.,
and let  $0<\dz<1$,   for any $\eta>0$. Then there exists a constant $C>0$ such
that
$$M^\sharp_{\vz,\dz,\eta}(Tf)(x)\le C
M_{\vz,\eta}(f)(x), \ \ {\rm a.e.}\ \ x\in \rz$$ for any smooth
 function $f$ with
compact support, where $\vz_\eta(Q)=(1+r)^\eta$ with $r=|Q|^{1/n}$, and
$$M_{\vz,\eta}f(x)=\dsup_{x\in Q}\dfrac 1{\vz_\eta(Q)|Q|}\dint_Q|f(y)|\,dy,$$
and $M_{\vz,\dz,\eta}^\sharp
f(x)=M^\sharp_{\vz,\eta}(|f|^\dz)^{1/\dz}(x),$ and the  sharp maximal operator $M^\sharp_{\vz,\eta}f(x)$ is defined by
$$\begin{array}{cl}
M_{\vz,\eta}^\sharp f(x)&:=\dsup_{x\in Q,r<1}\dfrac
1{|Q|}\dint_{Q_{x_0}}|f(y)-f_Q|\,dy+ \dsup_{x\in Q,r\ge 1}\dfrac
1{\vz_\eta(Q)|Q|}\dint_{Q_{x_0}}|f|\,dy\\
&\simeq\dsup_{x\in Q,r< 1}\dinf_C \frac
1{|Q|}\int_{Q_{x_0}}|f(y)-C|\,dy+ \dsup_{x\in Q,r\ge1}\dfrac 1{\vz_\eta |Q|}\dint_{Q_{x_0}}|f|\,dy
\end{array}$$
where $Q_{x_0}$ denotes cubes $Q(x_0,r)$ and $f_Q=\frac
1{|Q|}\int_Q f(x)dx$.
\end{lem}
\begin{lem}\label{l4.2.}\hspace{-0.1cm}{\rm\bf 4.2(\cite{t3}).}\quad
Let $0<p,\  \dz<\fz$ and $\wz\in A_\fz^{1,\fz}=A_\fz^{\rho,\fz}$ with $\rho=1$. There exists a
positive constant $C$ such that
$$\dint_\rz M_{\vz,\dz,\eta}
 f(x)^p\wz(x)dx\le C \dint_\rz M_{\vz,\dz,\eta}^\sharp
f(x)^p\wz(x)dx,$$
where $M_{\vz,\dz,\eta}
f(x)=M_{\vz,\eta}(|f|^\dz)^{1/\dz}(x),$
\end{lem}
From Lemmas 4.1 and 4.2, we have that for all $0<p<\fz$ and $\wz\in A^{1,\fz}_\fz$, for
any $\eta>0$, then there exists a constant $C$ depending only on
$\eta, p$ and the $A_\fz^{1,\fz}$ constant of $\wz$
such that
$$\|Tf\|_{L^p(\wz)}\le C\|M_{\vz,\eta}f\|_{L^p(\wz)}.$$
From this, we can get the vector-valued estimates (3.4) and (3.5) which  are new.

Next, we consider the  mulitilinear pseudo-differential operators, that is, $T$ is an $m$-linear operator such that $T$ are initially defined on the
$m$-fold product of Schwartz space ${\cal S}(\rz)$ and take their
values into the space of tempered distributions ${\cal S}'(\rz)$.
We will assume that the distributional kernel on $(\rz)^{m+1}$ of
the operator coincides away from the diagonal
$y_0=y_1=y_2=\cdots=y_m$ in $(\rz)^{m+1}$ with a function $K$ for
integer $m\ge 1$ so that
$$T(f_1,\cdots,f_m)(x)=\dint_{(\rz)^m}K(x,y_1,\cdots,y_m)f_1(y_1)
\cdots f_m(y_m)\,dy_1\cdots dy_m,$$ whenever $f_1,\cdots,f_m$ are
$C^\fz$ functions with compact support and
$x\not\in\bigcap_{j=1}^m supp f_j$. Moreover, we will assume that
the function $K$ satisfies the following estimates for any $N\ge 0$
$$|K(y_0,y_1,\cdots,y_m)|\le \dfrac
{C_N}{(1+\sum_{k,l=0}^m|y_k-y_l|)^{N}(\sum_{k,l=0}^m|y_k-y_l|)^{mn}},\eqno(4.8)$$ and, for some $\ez>0$ and any $N\ge 0$,
$$\begin{array}{cl}
|K(y_0,\cdots,y_j,\cdots, y_m)&-K(y_0,\cdots,y'_j,\cdots,y_l)|\\
&\le
\dfrac
{C_N|y_j-y'_j|^\ez}{(1+\sum_{k,l=0}^m|y_k-y_l|)^{N}(\sum_{k,l=0}^m|y_k-y_l|)^{mn+\ez}},\end{array}\eqno(4.9)$$
provided that $0\le j\le m$ and $|y_j-y'_j|\le \frac 12\ max_{\
0\le k\le m}|y_k-y'_k|$. When $N=0$ in (4.8) and (4.9), such kernels are called $m$-linear
Calder\'on-Zygmund kernels and the collections is denoted in
\cite{gt} by $m-CZK$. For these operators above, a boundedness
estimate
$$T:\ L^{q_1}\times\cdots\times L^{q_m}\to L^q,$$
for $1<q_1,\cdots,q_m<\fz$ and
$$\frac 1{q_1}+\cdots+\frac 1{q_m}=\frac 1q,\eqno(4.10)$$
implies the boundedness of the operator for all possible exponents
in such range of values.  Moreover, it will be important for
purpose the following end-point estimate also satisfied by such
operators:
$$\ L^{q_1}\times\cdots\times L^{q_m}\to L^{q,\fz},$$
for $1\le q_1,\cdots,q_m<\fz$ satisfying (4.10). In particular, it
will be relevant the case
$$\ L^{1}\times\cdots\times L^1\to
L^{1/m,\fz},$$ which extends the classical result in the linear
case  $T:\ L^1\to L^{1,\fz}$; see \cite {gt}.

A typical example, Let $T$ is a bilinear pseudo-differential operator with symbols belonging to $SB_{1,0}^0$; see \cite{b}. From \cite{x}, we know that the kernel of $T$ satisfies (4.8) and (4.9), and it is bounded from $L^1\to L^{1/2}\times L^{1/2}$; see \cite{b}.

We next give a  estimate for a mulitilinear pseudo-differential operator.

\begin{lem}\label{l4.3.}\hspace{-0.1cm}{\rm\bf 4.3.}\quad
Let $T$ be a mulitilinear pseudo-differential operator as above.
Let  $0<\dz<1/m$,   $\eta_j>0$ for $j=1,\cdots,m$, and
$\eta=\sum_{j=1}^m\eta_j$. Then there exists a constant $C>0$ such
that
$$M^\sharp_{\dz,\eta}(T\overrightarrow{f})(x)\le C\dprod_{j=1}^m
M_{\vz,\eta_j}(f_j)(x), \ \ {\rm a.e.}\ \ x\in \rz\eqno(4.11)$$ for any smooth
vector function $\overrightarrow{f}=(f_1,f_2,\cdots, f_m)$ with
compact support .
\end{lem}
{\it Proof.}\quad Let $\overrightarrow{f}$ be any smooth vector
function. Let $x\in Q=Q(x_0,r)$. Write each
$\overrightarrow{f}=\overrightarrow{f^0}+\overrightarrow{f^\fz}$,
where
$\overrightarrow{f^0}=\overrightarrow{f}\chi_{2Q}=(f_{1}\chi_{2Q},
\cdots,f_{m}\chi_{2Q})$. Set
$$C_Q=(T(\overrightarrow{f^\fz})_Q=T(f_1^\fz,\cdots,f_m^\fz)_Q
.$$ It is easy to see that
$$|T(\overrightarrow{f})(z)-C_Q|\le |T(\overrightarrow{f^\fz})(z)-C_Q|
+C_{qm}\dsum_m T(f_1^{r_1},\cdots,f_m^{r_m})(z),\eqno(4.12)$$ where
in the last sum each $r_j=0$ or $\fz$ and in each term there is at
least are $r_j=0$.

To prove (4.11), we consider two cases about $r$, that is, $r< 1$ and
$r\ge 1$.

Case 1. when  $r< 1$. Using the regularity of the kernel (4.9) and
 $0<\dz<1/m<1$, by Minkowski's inequality, we get
$$\begin{array}{cl}
\l(\dfrac
1{|Q|}\dint_Q\r.&\l.|T^\dz(\overrightarrow{f^\fz})(z)-C_Q^\dz|\,dz\r)^{1/\dz}\\
&\le\l(\dfrac
1{|Q|}\dint_Q|T(\overrightarrow{f^\fz})(z)-(T(\overrightarrow{f^\fz}))_Q
|^\dz\,dz\r)^{1/\dz}\\
&=\dfrac 1{|Q|}\dint_Q\l|\dfrac 1{|Q|}\dint_Q
T(\overrightarrow{f^\fz})(z)
-T(\overrightarrow{f^\fz})(y)\,dy\r|\,dz\\
&\le\dfrac {C}{|Q|}\dint_Q\dfrac
1{|Q|}\dint_Q\dint_{(\rz)^m\setminus (2Q)^m}\\
 &\quad\times
(K(z,\overrightarrow{w})-k(y,\overrightarrow{w}))
\dprod_{j=1}^mf_{j}(w_j)\,d\overrightarrow{w}\,dy|\,dz\\
&\le\dfrac {C_N}{|Q|}\dint_Q\dfrac
1{|Q|}\dint_Q\dint_{(\rz)^m\setminus (2Q)^m}\dfrac1{(1+|y-w_1|+\cdots+|y-w_m|)^N}\\
&\quad\times \dfrac{|y-z|^\ez} {(|y-w_1|+\cdots+|y-w_m|)^{mn+\ez}}
\prod_{j=1}^m|f_{j}(w_j)|\,d\overrightarrow{w}\,dydz\\
\end{array}$$
$$\begin{array}{cl}
&\le C_N |Q|^{\frac \ez n}\dprod_{j=1}^m\l(\dint_{\rz\setminus
2Q}\dfrac {|f_{j}(w_j)}{(1+|x-w_j|)^{\eta_j}|x-w_j|^{n+\frac\ez
m}}\,d\overrightarrow{w}\r)\\
&\le C_N\dprod_{j=1}^m|Q|^{\frac \ez {mn}}\dint_{\rz\setminus
2Q}\dfrac {|f_{j}(w_j)|}{(1+|x-w_j|)^{\eta_j}|x-w_j|^{n+\frac\ez m}}\,dw_j\\
&\le C_N\dprod_{j=1}^m M_{\vz,\eta_j}(f_j)(x),
\end{array}$$
if takeing $N=m\eta$.

The above computations gives the correct restimates for the first
term in the right hand side of (4.2). To estimate the sum in the
right hand side of (4.12) we distinguish between two kinds of term .
One, in which at least one of the $k_j=\fz$, and one final term in
which all the $k_j=0$. A typical representative of the first kind of
term is $T(f_1^\fz,\cdots,f_i^\fz,f^0_{i+1},\cdots,f_m^0)(z)$. Using
the notation $R_i=(\rz\setminus 2Q)^i\times (2Q)^{m-i}$, by
Minkowski's inequality, we have
$$\begin{array}{cl}
&\l(\dfrac
1{|Q|}\dint_Q|T(f_1^\fz,\cdots,f_i^\fz,f^0_{i+1},\cdots,f_m^0)(z)
|^\dz\,dz\r)^{1/\dz}\\
&\quad\le \dfrac
1{|Q|}\dint_Q|T(f_1^\fz,\cdots,f_i^\fz,f^0_{i+1},\cdots,f_m^0)(z)|\,dz\\
&\quad\le \dfrac
1{|Q|}\dint_Q\l|\dint_{(\rz)^m}k(x,\overrightarrow{y})
f_1^\fz(y_1),\cdots,f_i^\fz(y_i),f^0_{i+1}(y_{i+1}),\cdots,f_m^0(y_m)
\,d\overrightarrow{y}\r|\, dz\\
&\quad\le \dfrac C{|Q|}\dint_Q\dint_{R_i}\dfrac
{|f_{1}^\fz(y_1)\cdots f_{i,k}^\fz(y_i)f^0_{i+1}(y_{i+1})\cdots
f_{m}^0(y_m)|} {(1+|z-y_1|+\cdots+|z-y_m|)^{N}(|z-y_1|+\cdots+|z-y_m|)^{mn}}
\,d\overrightarrow{y}\, dz\\
&\quad\le \dfrac {C_N}{|Q|}\dint_Q\l(\prod_{j=l+1}^m\dint_{2Q}
|f_j(y_j)|\,dy_j\prod_{j=1}^i\dint_{\rz\setminus 2Q}\dfrac
{|f_{j}(y_j)|}{(1+|z-y_j|)^{\eta_j}|z-y_j|^{\frac {mn}i}}\,dy_j\r)\, dz\\
&\quad\le C_N\l(\dprod_{j=l+1}^m\dint_{2Q}
|f_j(y_j)|\,dy_j\dprod_{j=1}^i\dint_{\rz\setminus 2Q}\dfrac
{|f_{j}(y_j)|}{(1+|z-y_j|)^{\eta_j}|x-y_j|^{\frac {mn}i}}\,dy_j\r)\\
&\quad\le C_N\dprod_{j=l+1}^m\dint_{2Q} |f_j(y_j)|\,dy_j
\dprod_{j=1}^i\dint_{\rz\setminus 2Q}\dfrac
{|f_{j}(y_j)|}{(1+|z-y_j|)^{\eta_j}|x-y_j|^{\frac {mn}i}}\,dy_j\\
&\quad\le C_N\dprod_{j=l+1}^m M_{\vz,\eta_j}(f_j)(x)(|Q|^{\frac {m-i}i})^i
\dprod_{j=1}^i\dint_{\rz\setminus 2Q}\dfrac
{|f_{j}(y_j)|}{(1+|z-y_j|)^{\eta_j}|x-y_j|^{n+\frac {n(m-i)}i}}\,dy_j\\
&\quad\le C_N\dprod_{j=1}^m M_{\vz,\eta_j}(f_j)(x),
\end{array}$$
where we have used that $m>i$ and $N=m\eta$.

Applying Kolmogorov's estimate ( \cite {p}) to the term
$T(\overrightarrow{f^0})=T(f_1^0,\cdots,f_m^0)(z)$, we have
$$\begin{array}{cl}
\l(\dfrac
1{|Q|}\dint_Q|T(\overrightarrow{f^0})(z)|^\dz\,dz\r)^{1/\dz} &\le
C\|T(\overrightarrow{f^0})\|_{L^{1/m,\fz}}(Q,\frac
{dx}{|Q|})\\
&\le C\dprod_{j=1}^m\dfrac 1{|Q|}\dint_Q|f_j(z)|\,dz\le C\dprod_{j=1}^m M_{\vz,\eta_j}(f_j)(x),
\end{array}$$
since $T:\ L^1\times\cdots\times L^1\to L^{1/m,\fz}$; see \cite{gt}.

Case 2. When  $r\ge 1$. Similar to the proof of case 1. Taking
$N=m\eta$, then
$$\begin{array}{cl}
\l(\dfrac
1{\vz_\eta(Q)|Q|}\dint_Q|T^\dz(\overrightarrow{f^\fz})(z)|\,dz\r)^{1/\dz}
&\le\dfrac {C}{|Q|}\dint_Q\dint_{(\rz)^m\setminus
(2Q)^m}|K(z,\overrightarrow{w})|
\dprod_{j=1}^mf_{j}(w_j)\,d\overrightarrow{w}\,dy|\,dz\\
&\le C_N\dint_{(\rz)^m\setminus (2Q)^m}\dfrac{\prod_{j=1}^m|f_{j}(w_j)|\,d\overrightarrow{w}\,dy} {(|y-w_1|+\cdots+|y-w_m|)^{mn+N}}\\
&\le C_N \dprod_{j=1}^m\l(\dint_{\rz\setminus 2Q}\dfrac
{|f_{j}(w_j)}{|x-w_j|^{n+\frac N
m}}\,d\overrightarrow{w}\r)\\
&\le C_N\dprod_{j=1}^m\dint_{\rz\setminus
2Q}\dfrac {|f_{j}(w_j)|}{|x-w_j|^{n+\frac N m}}\,dw_j\\
&\le C_N\dprod_{j=1}^m M_{\vz,\eta_j}(|f_j|)(x).
\end{array}$$
The above computations gives the correct estimates for the first
term in the right hand side of (4.12). To estimate the sum in the
right hand side of (4.12) we distinguish between two kinds of term .
One, in which at least one of the $k_j=\fz$, and one final term in
which all the $k_j=0$. A typical representative of the first kind of
term is $T(f_1^\fz,\cdots,f_i^\fz,f^0_{i+1},\cdots,f_m^0)(z)$. Using
the notation $R_i=(\rz\setminus 2Q)^i\times (2Q)^{m-i}$, by
Minkowski's inequality, we have
$$\begin{array}{cl}
&\l(\dfrac
1{|Q|}\dint_Q|T(f_1^\fz,\cdots,f_i^\fz,f^0_{i+1},\cdots,f_m^0)(z)
|^\dz\,dz\r)^{1/\dz}\\
&\quad\le \dfrac
1{|Q|}\dint_Q|T(f_1^\fz,\cdots,f_i^\fz,f^0_{i+1},\cdots,f_m^0)(z)|\,dz\\
&\quad\le \dfrac
1{|Q|}\dint_Q\l|\dint_{(\rz)^m}k(x,\overrightarrow{y})
f_1^\fz(y_1),\cdots,f_i^\fz(y_i),f^0_{i+1}(y_{i+1}),\cdots,f_m^0(y_m)
\,d\overrightarrow{y}\r|\,dz\\
&\quad\le \dfrac {C_N}{|Q|}\dint_Q\l(\dint_{R_i}\dfrac
{|f_{1}^\fz(y_1)\cdots f_{i,k}^\fz(y_i)f^0_{i+1}(y_{i+1})\cdots
f_{m}^0(y_m)|}
{(1+(|z-y_1|+\cdots+|z-y_m|))^{N}(|z-y_1|+\cdots+|z-y_m|)^{mn}}
\,d\overrightarrow{y}\r)\, dz\\
&\quad\le \dfrac
{C_N}{|Q|}\dint_Q\l(\prod_{j=l+1}^m\vz(Q)^{-\eta_j}\dint_{2Q}
|f_j(y_j)|\,dy_j\prod_{j=1}^i\dint_{\rz\setminus 2Q}\dfrac
{|f_{j}(y_j)|}{|z-y_j|^{\frac {mn}i+\eta_j}}\,dy_j\r)\, dz\\
&\quad\le C_N\l(\dprod_{j=l+1}^m\dint_{2Q}
|f_j(y_j)|\,dy_j\dprod_{j=1}^i\dint_{\rz\setminus 2Q}\dfrac
{|f_{j}(y_j)|}{|x-y_j|^{\frac {mn}i+\eta_j}}\,dy_j\r)\\
&\quad\le C_N\dprod_{j=l+1}^m\dint_{2Q} |f_j(y_j)|\,dy_j
\dprod_{j=1}^i\dint_{\rz\setminus 2Q}\dfrac
{|f_{j}(y_j)|}{|x-y_j|^{\frac {mn}i+\eta_j}}\,dy_j\\
&\quad\le C_N\dprod_{j=l+1}^m M_{\vz,\eta_j}(f_j)(x)(|Q|^{\frac {m-i}i})^i
\dprod_{j=1}^i\dint_{\rz\setminus 2Q}\dfrac
{|f_{j}(y_j)|}{|x-y_j|^{n+\frac {n(m-i)}i+\eta_j}}\,dy_j\\
&\quad\le C_N\dprod_{j=1}^m M_{\vz,\eta_j}(f_j)(x),
\end{array}$$
where we have used that $m>i$ and $N=m\eta$.

Applying Kolmogorov's estimate ( \cite {p}) to the term
$T(\overrightarrow{f^0})=T(f_1^0,\cdots,f_m^0)(z)$, we have
$$\begin{array}{cl}
\l(\dfrac
1{\vz_\eta(Q)|Q|}\dint_Q|T(\overrightarrow{f^0})(z)|^\dz\,dz\r)^{1/\dz}
&\le
C\vz_\eta(Q)^{-1}\|T(\overrightarrow{f^0})\|_{L^{1/m,\fz}}(Q,\frac
{dx}{|Q|})\\
&\le C\dprod_{j=1}^m\dfrac 1{\vz_{\eta_j}(Q)|Q|}\dint_Q|f_j(z)|\,dz\\
&\le C\dprod_{j=1}^m M_{\vz,\eta_j}(f_j)(x),
\end{array}$$
since $T:\ L^1\times\cdots\times L^1\to L^{1/m,\fz}$; see \cite{gt}.

Hence, Lemma 4.3 is proved.\hfill$\Box$

Applying Lemma 4.2 and 4.3, we show that for $1<p<\fz$ and for all $\wz\in A_\fz^{1,\fz}$, for any $\eta_j>0$, $j=1,\cdots,m$,
$$\|T(f_1,\cdots,f_m)\|_{L^p(\wz)}\le C\|\dprod_{j=1}^m M_{\eta_j}(f_j)\|_{L^p(\wz)}.\eqno(4.13)$$
The scalar estimate (3.3) just (4.13). But the vector-valued inequalities (3.4) and (3.5) are new and immediately yield the following
result by applying H\"older's inequality and the norm inequalities for the maximal operator.
\begin{thm}\label{t4.4.}\hspace{-0.1cm}{\rm\bf 4.4.}\quad
Let $T$ be a mulitilinear pseudo-differential operator, $1\le p_1,\cdots,$ $ p_m<\fz$, $1<q_1,\cdots, q_m<\fz$ and $0<p,q<\fz$ such that
$$\frac 1p=\frac 1{p_1}+\cdots+\frac 1{p_m},\quad  \frac 1q=\frac 1{q_1}+\cdots+\frac 1{q_m}.$$
If $1<p_1,\cdots, p_m<\fz$ and $\wz\in A_{p_1}^{1,\fz}\bigcap\cdots\bigcap A_{p_m}^{1,\fz}$, then
$$\||T\overrightarrow{f}|_q\|_{L^p(\wz)}\le C\dprod_{j=1}^m\||f|_{q_j}\|_{L^{p_j}(\wz)}.\eqno(4.14)$$
If at leat one $p_j=1$ and $\wz\in A_1^{1,\fz}$, then
$$\||T\overrightarrow{f}|_q\|_{L^{p,\fz}(\wz)}\le C\dprod_{j=1}^m\||f|_{q_j}\|_{L^{p_j}(\wz)}.\eqno(4.15)$$
Moreover, inequalities (4.14) and (4.15) hold with $T^*$ in place $T$, where $T^*$ is the dual operator of $T$.

\end{thm}

{\bf Remark.} We will continue to study weighted inequalities for mulitilinear pseudo-differential operators in the forthcoming paper.

 LMAM, School of Mathematical  Science

 Peking University

 Beijing, 100871

 P. R. China

\bigskip

 E-mail address:  tanglin@math.pku.edu.cn


\begin{center} {\bf References}\end{center}
\begin{enumerate}
\vspace{-0.3cm}
\bibitem[1]{aj} K. Andersen and R. John,
Weighted inequalities for vector-valued maximal functions and singular integrals, Studia. math. T. LXIX. (1980), 19-31.
\vspace{-0.3cm}
\bibitem[2]{bp} A. Benedek and R. Panzone,
The space $L^p$ with mixed norm, Duke Math. J. 28(1961), 301-324.
\vspace{-0.3cm}
\bibitem[3]{b}
 $\acute{A}$. B$\acute{e}$nyi and R. Torres,  Symbolic calculus and the transposes of bilinear pseudodifferential
operators, Comm. Par. Diff. Eq. 28 (2003), 1161-1181.
\vspace{-0.3cm}
\bibitem[4]{bhs1} B. Bongioanni, E. Harboure and O. Salinas,
Commutators of Riesz transforms related to Schr\"odinger operators, J. Fourier Ana Appl. 17(2011), 115-134.
\vspace{-0.3cm}
\bibitem[5]{bhs2} B. Bongioanni, E. Harboure and O. Salinas,
Class of weights related to Schr\"odinger operators, J. Math. Anal. Appl. 373(2011), 563-579.
 \vspace{-0.3cm}
\bibitem[6]{cmp} D. Cruz-Uribe, J. M. Martell and C. P\'erez,
Extrapolation from $A_\fz$ weights and applications, J. Funct. Anal. 213(2004), 412-439.
 \vspace{-0.3cm}
\bibitem[7]{dz}J. Dziuba\'{n}ski and J. Zienkiewicz,
 Hardy space $H^1$ associated to Schr\"{o}dinger operator with potential
satisfying reverse H\"{o}lder inequality, Rev. Math. Iber. 15
(1999), 279-296. \vspace{-0.3cm}
\bibitem[8]{dz1}J. Dziuba\'{n}ski, G. Garrig\'{o}s, J. Torrea and J.
Zienkiewicz, $BMO$ spaces related to Schr\"{o}dinger operators with
potentials satisfying a reverse H\"{o}lder inequality, Math. Z.
249(2005), 249 - 356. \vspace{-0.3cm}
\bibitem[9]{d} J. Duoandikoetxea,
Fourier Analysis, in: Graduate Studies in Mathematics, Vol. 29,
American Mathematical Society, Providence, RI, 2000.
 \vspace{-0.3cm}
\bibitem[10]{fs}C. Fefferman and E. Stein,
Some maximal inequalities, Amer. J. Math. 93(1971), 107-115.
\vspace{-0.3cm}
\bibitem[11]{gt} L. Grafakos and R. Torres,
Multilinear Calder\'on-Zygmund theory, Adv. Math. 165(2002),
124-164.
\vspace{-0.3cm}
\bibitem[12]{glp}
Z. Guo, P. Li and L. Peng, $ L^p$ boundedness of commutators of
Riesz transforms associated to  Schr\"{o}dinger operator, J. Math.
Anal and Appl. 341(2008), 421-432.
\vspace{-0.3cm}
\bibitem[13]{g} J. Garc\'ia-Cuerva,
An extrapolation theorem in the theory of $A_p$-weights,
Proc. Amer. Math. Soc. 87(1983), 422-426.
 \vspace{-0.3cm}
\bibitem[14]{gr} J. Garc\'ia-Cuerva and J. Rubio de Francia,
Weighted norm inequalities and related topics, Amsterdam- New
York, North-Holland, 1985.
\vspace{-0.3cm}
\bibitem[15]{j} P. Jones,
Factorization of $A_p$ weights, Ann of Math. 111(1980), 511-530.
 \vspace{-0.3cm}
\bibitem[16]{m} B. Muckenhoupt,
Weighted norm inequalities for the Hardy maximal functions, Trans.
Amer. Math. Soc. 165(1972), 207-226.
 \vspace{-0.3cm}
\bibitem[17]{p}C. P\'erez,
Endpoint estimates for commutators of singular integral operators,
J. Funct. Anal. 128(1995), 163-185. \vspace{-0.3cm}
\bibitem[18]{rr} M. M. Rao and Z. D. Ren,
Theory of Orlicz spaces, Monogr. Textbooks Pure Appl. Math.146,
Marcel Dekker, Inc., New York, 1991.
\vspace{-0.3cm}
\bibitem[19]{ru}J. Rubio de Francia,
Factorization theory and $A_p$ weights, Amer. J. Math. 106 (1984),
533-547. \vspace{-0.3cm}
\bibitem[20]{s1}Z. Shen,
$L^p$ estimates for Schr\"odinger operators with certain
potentials, Ann. Inst. Fourier. Grenoble, 45(1995), 513-546.
\vspace{-0.3cm}
\bibitem[21]{s}  E. M. Stein,
Harmonic Analysis: Real-variable Methods, Orthogonality, and
Oscillatory integrals. Princeton Univ Press. Princeton, N. J.
1993.
\vspace{-0.3cm}
\bibitem[22]{t} M. Taylor,
Pseudodifferential operators and nonlinear PDE. Boston:
Birkhau-ser, 1991.
\vspace{-0.3cm}
\bibitem[23]{t1} L. Tang,
Weighted norm inequalities for  Schr\"odinger type operators, preprint.
\vspace{-0.3cm}
\bibitem[24]{t2} L. Tang,
Weighted norm inequalities for commutators of Littlewood-Paley functions related to Schr\"odinger operators, preprint.
\vspace{-0.3cm}
\bibitem[25]{t3} L. Tang,
Weighted norm inequalities for pseudo-differential
operators with smooth symbols and their commutators, preprint.
\vspace{-0.3cm}
\bibitem[26]{w} D. Watson,
Vector-valued inequalies, factorization, and extrapolation for a family of rough operators, J. Funct. Anal. 121(1994), 389-415.
\vspace{-0.3cm}
\bibitem[27]{x} J. Xiao, Y. Jiang and W. Gao,
Bilinear Pseudo-differential Operators on Local
Hardy Spaces, to appear in Acta Math. Sin. (Engl. Ser.).
\vspace{-0.3cm}
\bibitem[28]{z}  J. Zhong,
Harmonic analysis for some Schr\"odinger type operators, Ph.D. Thesis. Princeton University,
  1993.

\end{enumerate}
\end{document}